\documentclass[11pt]{article}

\usepackage{graphicx}
\usepackage{amssymb}
\usepackage{enumerate}

\newcommand{\Ceq}{\stackrel{+}{=}}

\newcommand\independent{\protect\mathpalette{\protect\independenT}{\perp}}
\def\independenT#1#2{\mathrel{\rlap{$#1#2$}\mkern2mu{#1#2}}}

\newcommand{\cA}{{\cal A}}

\newcommand{\R}{{\mathbb R}}

\newcommand{\N}{{\mathbb N}}
\newcommand{\Z}{{\mathbb Z}}

\newcommand{\beq}{\begin{equation}}
\newcommand{\eeq}{\end{equation}}
\newcommand{\beqy}{\begin{eqnarray}}
\newcommand{\eeqy}{\end{eqnarray}}

\newtheorem{Postulate}{Postulate}
\newtheorem{Definition}{Definition}
\newtheorem{Lemma}{Lemma}
\newtheorem{Theorem}{Theorem}

\newenvironment{Definition*}{{\bf Definition}}{}

\date{April 23, 2008}

\title{Causal inference using \\the algorithmic Markov condition}

\author{Dominik
Janzing \,and Bernhard Sch\"olkopf\thanks{{\tt email: \{dominik.janzing,bernhard.schoelkopf\}@tuebingen.mpg.de}}
\\
 \small
Max Planck Institute for Biological Cybernetics,\\
\small Spemannstr. 38,\\
72076 T\"ubingen, Germany
}

\begin{document}
\maketitle

\abstract{Inferring the causal structure that links $n$ observables is usually based
upon detecting statistical dependences and choosing simple graphs that
make the joint measure Markovian. Here we argue why causal inference is also possible when only
single observations are present.

We develop a theory how to generate causal graphs explaining
similarities between single objects. To this end, we replace the notion of  
conditional stochastic independence in the causal Markov condition 
with the vanishing of conditional {\it algorithmic} mutual information 
and describe the corresponding causal inference rules.

We explain why a consistent reformulation of causal inference in terms of 
algorithmic complexity implies a new inference principle that takes into account also the
complexity of conditional probability densities, making it possible
to select among Markov equivalent causal graphs.  This insight provides  a theoretical foundation of
a heuristic principle proposed in earlier work.

We also discuss how to replace Kolmogorov complexity with {\it decidable} 
complexity criteria. This can be seen as an algorithmic analog of 
replacing the empirically  undecidable question of statistical independence with
practical independence tests that are based on implicit or  explicit assumptions
on the underlying distribution.
}

\newpage

\tableofcontents

\section{Introduction to causal  inference from statistical data}

\label{Sec:I}

Causal inference from statistical data has attracted increasing interest in the past decade. In  contrast to traditional statistics where 
statistical dependences are only taken to prove that some kind of relation between random variables exists, 
causal inference methods in machine learning are explicitly designed to generate hypotheses on causal  directions automatically based upon
statistical independence tests \cite{Pearl:00,Spirtes}. The crucial assumption connecting statistics with causality is the causal Markov condition explained below  after we have introduced some  notations and terminology.

We denote  random variables by capitals and their values by the corresponding lowercase letters.
Let $X_1,\dots,X_n$ be random variables and $G$ be a directed acyclic graph (DAG) representing the causal structure where an arrow from node $X_i$ to node $X_j$
indicates a direct causal effect. 
Here the term {\it direct} is understood with respect to the chosen set  of  variables in 
the sense that the information flow between the two variables considered is not performed via using one or more of the other variables
as intermediate  nodes. 
We will next briefly rephrase the postulates that are required in the statistical theory of inferred causation \cite{Spirtes,Pearl:00}.

\subsection{Causal Markov condition}

\label{Subsec:MC}

When we  consider  the causal structure that links $n$ random variables
 ${\cal V}:=\{X_1,\dots,X_n\}$
we will implicitly  assume that ${\cal V}$ is causally sufficient in the sense that all common  causes of two  variables
in ${\cal V}$ are also in ${\cal V}$. 
Then
a causal hypothesis $G$ is only acceptable as potential causal structure if the joint distribution $P$ of
$X_1,\dots,X_n$ satisfies the Markov condition with respect to $G$. There  are several formulations of the Markov condition that
are  known to coincide under some technical  condition (see Lemma~\ref{equivMKSt}). 
We will first introduce the following version which is sometimes referred to as the {\it parental} or the {\it local} Markov condition  
\cite{LauritzenDawid}.  

To  this end,  we introduce the following notations. $PA_j$  is the set of  parents  of $X_j$ 
and $ND_j$ the set  of non-descendants of $X_j$ except itself. If $S,T,R$  are sets of  random  variables,
$S \independent T\, |R$ means $S$ is statistically independent of  $T$, given $R$.

\begin{Postulate}[statistical  causal Markov condition, local]${}$\\ \label{PMC}
If a directed acyclic graph $G$ formalizes the causal structure among the random variables $X_1,\dots, X_n$. Then 
\[
X_j \independent ND_j  \,|PA_j\,,
\]
for all $j=1,\dots,n$.  
\end{Postulate}

We call this postulate the {\it statistical} causal Markov condition  because we  will  later introduce
an algorithmic version. 
The fact that conditional irrelevance not  only occurs in  the  context of {\it statistical} dependences
has been emphasized in the literature (e.g. \cite{Lauritzen,Pearl:00}) in the context of describing
abstract properties (like semi-graphoid axioms) of
the relation $  \cdot \independent \cdot \,|\cdot$. We will therefore state  the causal Markov condition also
in an abstract form that does not refer to any specific notion of conditional informational irrelevance:

\begin{Postulate}[abstract causal  Markov condition, local]${}$\\ \label{abstractMC} 
Given all the direct causes of an observable $O$, its non-effects provide no additional information on $O$. 
\end{Postulate} 

Here,
observables denote something in the real world that can be observed and the observation of which can be 
formalized in terms of a mathematical language. 
In this paper,
observables will either be random variables (formalizing statistical quantities) or they will be 
strings (formalizing the description of objects).
Accordingly, information will be {\it statistical} or {\it algorithmic}  mutual information, respectively.  

The importance of  the causal Markov  condition lies in the fact that it links causal terms like ``direct causes'' and 
``non-effects'' to informational relevance of observables.
The local Markov condition  is rather intuitive because it
echoes the fact that the information flows from direct causes to their effect
and every dependence between a node and its non-descendants involves the direct
causes. However, the independences postulated  by the local Markov condition imply  additional
independences. It is therefore hard to decide whether
an independence must hold for a Markovian distribution or not, solely on the basis  of the local formulation.  
In contrast, the global Markov condition makes the complete set of independences obvious. 
To state it we first have to introduce the following graph-theoretical  concept.

\begin{Definition}[d-separation]${}$\\ \label{dsep} 
A path  $p$ in  a DAG is said to be d-separated (or blocked) by a set of nodes $Z$ if and only 
if 
\begin{enumerate}
\item $p$ contains a chain $i\rightarrow m \rightarrow j$ or fork $i\leftarrow m  \rightarrow j$ such that the middle node $m$ is in $Z$, or
\item  $p$  contains an inverted fork (or collider) $i\rightarrow m \leftarrow j$ such that the middle node $m$ is not in $Z$ and such that
 no descendant of $m$ is in  $Z$.
\end{enumerate}
A set $Z$ is said to d-separate $X$ from $Y$ if and only if  $Z$ blocks every (possibly undirected) 
path from a node in $X$ to a node in $Y$.
\end{Definition}

The following Lemma shows that d-separation is  the correct condition for deciding 
whether an independence is implied by the local Markov condition \cite{Lauritzen}, Theorem~3.27.

\begin{Lemma}[equivalent Markov conditions]${}$\\ \label{equivMKSt}
Let 
$P(X_1,\dots,X_n)$ have a density $P(x_1,\dots,x_n)$  with respect to a product measure. 
Then the following three statements are equivalent:

\begin{enumerate}[{\em I.}]

\item {\bf Recursive form:} $P$ admits the factorization
\begin{equation}\label{RecSt}
P(x_1,\dots,x_n)=\prod_{j=1}^n P(x_j|pa_j)\,,
\end{equation}
where $P(.|pa_j)$ is shorthand for the conditional probability density, given 
the values of all parents of  $X_j$.

\item {\bf Local (or parental) Markov condition:} 
for every node  $X_j$  we have
\[
X_j \independent ND_j \,| PA_j\,,
\]
i.e.,
it is conditionally independent of its non-descendants (except itself), given its parents.

\item {\bf Global Markov condition:} 
\[ 
S \independent T \,| R
\]
 for all three sets $S,T,R$ of nodes
for which $S$ and $T$ are d-separated by $R$.

\end{enumerate}

Moreover, 
the local and the global Markov condition are equivalent even if $P$ does not have a density with respect to a product measure. 
\end{Lemma}

The  conditional  densities $P(x_j|pa_j)$ are also called the  {\it Markov kernels} relative to the hypothetical causal  graph $G$. 
It is important to note that every choice of Markov kernels define a Markovian density $P$, 
 i.e., the Markov kernels define exactly the set of free parameters remaining after the causal structure has  been specified. 

To select graphs among all those that render $P$ Markovian, we also need an additional postulate:

\begin{Postulate}[causal faithfulness]${}$\\ \label{faithfulness}
Among all graphs $G$ for  which $P$ is Markovian, prefer the ones for which all
the observed conditional independences in the joint measure $P(X_1,\dots,X_n)$  are imposed
by the Markov condition.
\end{Postulate}

The idea 
is that the set of observed independences is typical  for the causal structure under consideration rather than being the result  of specific choices of the Markov kernels. This becomes even more intuitive  when we restrict our attention to
random variables with finite value set and observe that the values $P(x_j|pa_j)$ then define a natural parameterization of the 
set of Markovian distributions in a finite dimensional space. The non-faithful distributions form a 
submanifold of lower dimension, i.e., a
set of Lebesgue measure zero \cite{Meek}. 
They therefore almost surely don't occur if we assume that ``nature chooses'' the Markov  kernels for the different nodes independently according to
some density on the parameter space. 

The above ``zero Lebesgue measure argument''   is close to the spirit of  Bayesian approaches  \cite{Heckerman1999a}, where priors on the set  of  Markov kernels
are  specified for every possible hypothetical causal DAG and causal inference is performed by maximizing
posterior probabilities for hypothetical DAGs, given the  observed data.
This procedure leads to an
{\it implicit} preference of faithful structures in the  infinite sampling limit given some natural conditions for the  priors
on the parameter space. The assumption that ``nature chooses Markov kernels
independently'', which is also part of the Bayesian approach, will turn out to be closely related to the algorithmic Markov condition postulated  in this paper.

We now discuss  the justification  of the 
statistical causal Markov condition because we will later justify
the algorithmic Markov condition in a similar way. 
To this end, we introduce functional models \cite{Pearl:00}:

\begin{Postulate}[functional model of causality]${}$\\  \label{functional}
If a directed acyclic graph $G$  formalizes the causal relation between the random variables $X_1,\dots,X_N$ 
then 
every $X_j$ can be written as  a deterministic function of $PA_j$ and a noise variable $N_j$\,,
\[
X_j=f_j(PA_j,N_j)\,,
\]
where all $N_j$ are jointly independent.
\end{Postulate}

\noindent
Then we have \cite{Pearl:00}, Theorem~1.4.1:

\begin{Lemma}[Markov condition in functional models]${}$\\ \label{funcM} 
Every joint distribution $P(X_1,\dots,X_n)$ generated according to 
the functional model in Postulate~\ref{functional}
satisfies the local and the global Markov condition relative to $G$.
\end{Lemma}

We rephrase the proof in \cite{Pearl:00} because our proof for the algorithmic version will rely on the  same  idea.

\vspace{0.3cm}
\noindent
Proof of Lemma~\ref{funcM}: extend $G$  to a graph $\tilde{G}$  with nodes $X_1,\dots,X_n,N_1,\dots,N_n$ that additionally  contains
an arrow from each $N_j$ to $X_j$.  The given joint distribution of noise variables
induces a joint distribution
\[
\tilde{P}(X_1,\dots,X_n,N_1,\dots,N_n)\,,
\]
that
satisfies the local Markov  condition with respect  to $\tilde{G}$: first, every $X_j$ is completely determined by its parents
making the condition trivial.
Second, every $N_j$ is parentless and thus we have to check that it is (unconditionally) independent of its non-descendants. The latter
are deterministic functions of $\{N_1,\dots,N_n\} \setminus \{N_j\}$. Hence the independence follows from the joint independence of
all $N_i$. 

By Lemma~\ref{equivMKSt}, $\tilde{P}$ is also globally Markovian w.r.t. $\tilde{G}$. 
Then we observe that $ND_j$ and $X_j$ are d-separated in $\tilde{G}$ (where the parents and non-descendants are defined 
with respect  to  $G$). Hence $P$ satisfies the local Markov  condition w.r.t. $G$ and hence also the global Markov condition.
$\Box$

\vspace{0.3cm}
Functional models formalize the idea that the outcome of an experiment is completely
determined by the values of all relevant parameters where the only uncertainty stems from the fact that 
some of these parameters are hidden.
Even though this kind of determinism is in contrast with the commonly accepted
interpretation of  quantum mechanics \cite{Omnes}, we still consider functional models 
as a helpful framework for discussing causality in real life since quantum mechanical laws refer mainly to phenomena in micro-physics.

Causal inference using the Markov condition and the faithfulness assumption 
has been implemented as causal learning algorithms \cite{Spirtes}. 
The following fundamental limitations of these methods
deserve our further attention:

\begin{enumerate}

\item {\it Markov equivalence:} There are only few cases where the inference rules provide unique  causal graphs. Often 
one ends up with a large
class of Markov equivalent  graphs, i.e., graphs that entail the same set of independences. For this reason, additional inference rules
are desirable. 

\item {\it Dependence on i.i.d.~sampling:} the whole setting of causal inference relies on
the  ability to sample repeatedly and independently from the same joint distribution $P(X_1,\dots,X_n)$. As opposed to this assumption, 
causal inference in real life also deals  with probability distributions that change in time and often one infers 
causal relations among single observations without referring to statistics at all.

\end{enumerate}

The idea of this paper is to  develop a theory of probability-free causal inference
that helps to construct causal hypotheses based on similarities of {\it single} objects.
Here, similarities will be defined by 
comparing the length of  the
shortest description of single  objects to the  length of their shortest joint description. 
Despite the analogy to 
causal inference from statistical data (which is due to known analogies between statistical and algorithmic information theory)
our theory also implies new {\it statistical} inference rules. 
In other words, our approach to address weakness  2 also yields new methods to address 1.

The paper is structured as  follows. In the remaining part of this Section, i.e., Subsection~\ref{Subsec:plMK}, 
we describe recent approaches from the  literature to causal inference from statistical
data that address problem 1 above. In Section~\ref{Sec:Ind} we develop the general theory on inferring causal relations 
among individual objects based on algorithmic information. This framework appears, at first sight,  as a straightforward adaption
of the statistical framework (using well-known correspondences between statistical and algorithmic information theory).
However,
 Section~\ref{Sec:NewInf} describes that this implies novel causal inference rules for {\it statistical} 
inference because {\it non-statistical} algorithmic dependences can even occur in data that were obtained from statistical sampling.
In Section~\ref{Sec:Dec} we describe how to replace causal inference rules based on the uncomputable {\it algorithmic information}
with decidable criteria that are still motivated  by the uncomputable idealization. 

The table in fig.~\ref{analog}  summarizes the analogies
between the theory of statistical and the theory of algorithmic causal inference described in this paper. 
The differences, however, which are the main subject of
Sections~\ref{Sec:NewInf} to \ref{Sec:Dec}, can hardly be represented in the table.

\vspace{0.3cm}
\begin{figure}
\begin{tabular}{clcccc}
&  && {\bf statistical} &&  {\bf algorithmic} \\
&  &&     &&  \\
&observables & & random variables && sequences of strings \\
& (vertices  of a DAG)  &&    && \\
&   && && \\
&observations && i.i.d.~sampled data && strings  \\  
&   &&    && \\
&conditional independence& \quad\quad&  $X \independent Y\,|Z$ && $x\independent y\,|z$ \\
&          &&  $\stackrel{..}{\Updownarrow}$ &&  $\stackrel{..}{\Updownarrow}$ \\
&    &&     $I(X;Y|Z)=0$     && $I(x:y|z) \Ceq 0$ \\  
&    && && \\
&I. recursion formula && $P(x_1,\dots,x_n)$ && $K(x_1,\dots,x_n)$ \\
&                  && $=$ && $=$ \\
& && $\prod_j p(x_j|pa_j)$ && $\sum_j K(x_j|pa_j^*)$ \\
& && &&  \\
&II. local Markov condition&& $X_j \independent  ND_j\,|PA_j$  && $x_j\independent nd_j  \,|pa_j^*$  \\ 
& &  & \\
&III. global Markov&&  d-separation && d-separation  \\
&condition  &&  $\Rightarrow$  && $\Rightarrow$ \\
&  &&   statistical independence &&  algorithmic independence \\
& && && \\
&equivalence of I-III&& Theorem~3.27  &&Theorem~\ref{equiAMK} \\
& && in \cite{Lauritzen} &&   \\
& &&    &&  \\
&functional models&& Section 1.4  && Postulate~\ref{algoFunc} \\
&  &&  in \cite{Pearl:00}     && \\
&  & &  \\
&functional models&&  Theorem~1.4.1  && Theorem~\ref{algoImplM} \\
&imply Markov condition&&   in  \cite{Pearl:00} && \\
& && && \\
& decidable dependence &&  assumptions on    &&  Section~\ref{Sec:Dec}    \\
&criteria &&        joint distribution          &&   \\
& && &&   
\end{tabular} 
\caption{\label{analog}{\small Analogies between statistical and algorithmic causal inference}} 
\end{figure}

\subsection{Seeking for new statistical inference rules}

\label{Subsec:plMK}

In  \cite{SunLauderdale} and \cite{NeuroComp}  we have proposed causal  inference rules that are based on the idea  that
the factorization of $P({\rm cause},{\rm effect})$  into $P({\rm effect}|{\rm cause})$ and $P({\rm cause})$
typically leads to simpler terms than the ``artificial'' factorization into $P({\rm effect})P({\rm cause}|{\rm effect})$. 
The generalization  of this principle reads:
Among all graphs $G$ that render $P$ Markovian
prefer the one for which the decomposition  in eq.~(\ref{RecSt})
yields the simplest Markov kernels.  
We have
called this vague idea the ``principle of plausible Markov kernels''. 

\begin{figure}
\centerline{ \includegraphics[scale=0.25]{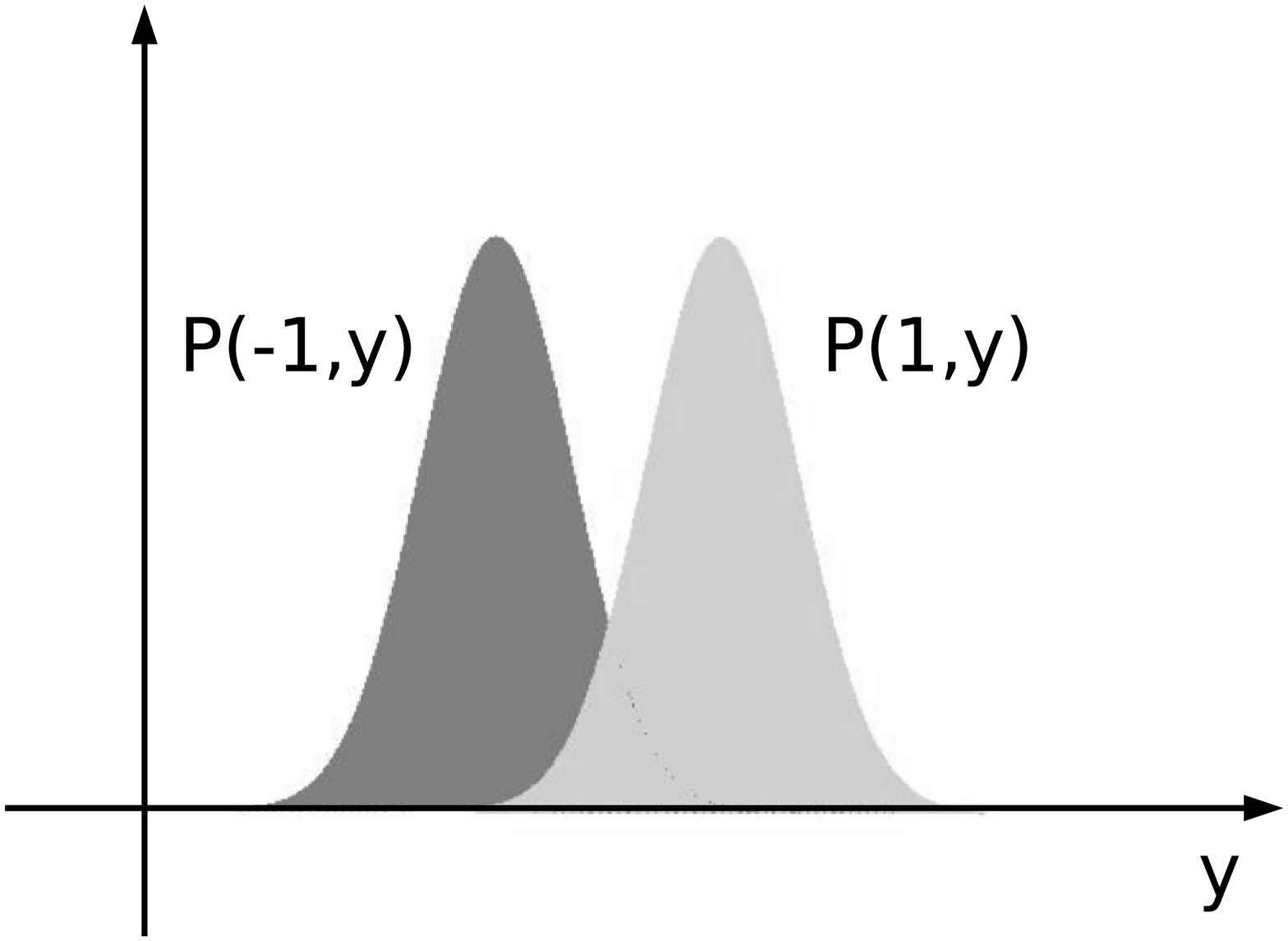}}
\caption{{\small Observed joint distribution of $X$ and $Y$ consisting  of  two Gaussians of equal width
shifted against each other.}\label{GaussianOrig}} 
\end{figure}

Before we describe several options to define simplicity we 
describe a simple example to illustrate the idea. 
Assume we have observed that a binary variable $X$ (with  values $x=-1,1$)
and a continuous variable $Y$ (with  values in  $\R$)
are distributed according to a mixture of two Gaussians (see fig.~\ref{GaussianOrig}). Since this will simplify  the
further discussion let us assume that the two  components are equally weighted, i.e.,
\[
P(x,y)= \frac{1}{2} \frac{1}{\sigma \sqrt{2\pi}}e^{-\frac{(y-\mu-x \lambda)^2}{2\sigma^2}} \,,
\]
where $\lambda$ determines the shift of the mean caused by switching between $x=1$ and $x=-1$.

The marginal $P(Y)$  is given by
\begin{equation}\label{marg}
P(y)=\frac{1}{2}\frac{1}{\sigma \sqrt{2\pi}}\left( e^{-\frac{(y-\mu+\lambda)^2}{2\sigma^2}} +  e^{-\frac{(y-\mu-\lambda)^2}{2\sigma^2}}\right)\,.
\end{equation}
One will  prefer  the causal structure $X\rightarrow Y$ compared  to $Y\rightarrow X$  because
the former explains in a natural  way why $P(Y)$ is bimodal: the effect of $X$ on $Y$ is simply to shift the Gaussian
distribution by $2 \lambda$. In the latter model the bimodality of $P(Y)$ remains unexplained. 
To prefer one causal model to  another one because the corresponding conditionals are simpler seems to be a natural  
application of Occam's Razor. However,  
Section~\ref{Sec:NewInf} will show that such an inference rule also follows from the theory developed in the present paper 
when simplicity is meant
in the sense of low Kolmogorov complexity. In the remaining part of this section we will sketch some approaches
to implement the ``principle of plausible Markov kernels'' in practical applications. 

In \cite{SunLauderdale} we have  defined a family  of ``plausible  Markov kernels''
by  conditionals $P(X_j|PA_j)$ that  are second order exponential models,  i.e.,
$\log P(x_j|pa_j)$ is a polynomial of order two in the variables $\{X_j\}\cup \{PA_j\}$ up to some
additive partition function  (for normalization) that depends only on the variables $PA_j$.
For every hypothetical causal graph, one thus obtains a  family of  
``plausible joint distributions $P(X_1,\dots,X_n)$'' that are  products of the plausible Markov kernels.
Then we prefer the causal direction for which the plausible joint distributions provide the best fit for the
given observations. 

In \cite{NeuroComp} we have proposed the following principle for causal inference:
Given a joint distribution of  the random variables $X_1,\dots,X_n$,  prefer
a causal structure for which 
\begin{equation}\label{Csum}
\sum_{j=1}^n C(P(X_j|PA_j))
\end{equation}
is minimal, where $C$ is some complexity measure on conditional probability densities.

There is also another recent proposal for new inference rules that refers to a related simplicity assumption,
though formally quite different from the ones above. The authors of \cite{Kano2003} observe that there are joint  distributions of $X_1,\dots,X_n$ 
that can  be explained by a linear model with additive non-Gaussian noise for one causal  direction but require non-linear 
causal influence for the other causal  directions. For real data they prefer the causal graph for which the observations are closer 
to the linear model. 

To justify the belief that conditionals that correspond to the true causal direction tend to be simpler than non-causal 
conditionals (which is common to all the  approaches above) is one of the main goals of  this paper.

\section{Inferring causal relations among individual objects}

\label{Sec:Ind}

It has  been emphasized \cite{Pearl:00} that the application of causal inference principles 
often benefits from the non-determinism of causal  relations between the observed  random variables.
In contrast, human learning in real-life often is about quite deterministic relations. Apart from that, the most important difference
between human causal learning and  the inference rules in \cite{Spirtes,Pearl:00} is that the former is also about
causal relations among {\it single} objects and does not  necessarily require sampling. 
Assume, for instance, that the comparison of two texts show similarities (see e.g. \cite{BennettScience}) such that the author of the text that appeared later is blamed to have
copied it from the other one or both are blamed to have copied from a third one. 
The statement that the texts are similar could be based on a statistical analysis of
the occurrences of certain words or letter sequences. However, such kind of simple statistical tests can fail in  both directions:
In Subsection~\ref{Subsec:Ind} (before Theorem~\ref{equiAMK}) 
we will discuss an example showing that they can erroneously infer causal relations even though they do not exist. This is because
parts that are common  two both objects,  e.g., the two texts, are only suitable to prove 
a causal link if they are not ``too straightforward'' 
to come up with.

On the other  hand, causal relations can  generate 
similarities between texts
for   which every {\it efficient} statistical analysis  
is believed to fail. We will describe an idea from cryptography to show this.
A  cryptosystem is called ROR-CCA-secure (Real or Random under Chosen Ciphertext Attacks) 
if there is no efficient method to decide whether 
a text is random or the encrypted version of some {\it known} text without knowing the key \cite{Hofheinz}.
Given that there are ROR-CCA-secure schemes (which is unknown but believed by cryptographers) 
we have a causal relation leading to similarities that are not detected by any kind of simple counting statistics.
However,  once an attacker has found the key (maybe  by exhaustive search), he recognizes similarities between 
the encrypted text and the plain text and infers a causal  relation. 
This already suggests two things: (1) detecting similarities involves {\it searching} over 
potential rules how properties of one object can be algorithmically  derived from the structure of the other. 
(2) It is likely that
inferring causal relations therefore relies on {\it computationally infeasible}  
decisions (if computable at all) on whether two objects 
have information in common or not.

\subsection{Algorithmic mutual information}

\label{Subsec:AM}

We will now describe how the information one object provides about the other can  be  measured in terms  of
Kolmogorov complexity.
We start with some notation and terminology. Below, strings will always be binary strings since
every description given in terms of a different  alphabet can be converted into a binary word. The set of binary strings
of arbitrary length will be denoted by $\{0,1\}^*$.
Recall that the Kolmogorov complexity
 $K(s)$ of a string $s\in \{0,1\}^*$ is defined as the length of the shortest program 
that generates $s$ using a previously defined universal Turing machine 
\cite{Solomonoff,Solomonoff2,KolmoOr,Chaitin,ChaitinF,Cover,Vitanyi97}.
The conditional Kolmogorov complexity $K(t|s)$  \cite{Cover} 
of a string $t$ given another string $s$ is the length of the shortest program that can generate
$t$ from $s$. 
In order to keep our notation simple  we use $K(x,y)$ to refer to the complexity of  the concatenation of $x,y$.  

We will mostly have equations that are valid only 
up to additive constant terms in the sense that the difference between both sides does  not depend on the strings
involved in the equation (but it may depend on the Turing machines they refer to). To indicate  such constants  we denote  the
corresponding  equality by $\stackrel{+}{=}$ and  likewise for inequalities. In this context it  is important to note  
that the number $n$ of nodes of the causal graph is  considered to be a constant. 
Moreover, for every string $s$ we define $s^*$ as its shortest description. If the latter is not unique, we consider  the first one in an  lexicographic order. 
It is necessary to distinguish between  $K(\cdot |s)$ and $K(\cdot| s^*)$. This is because
there is a trivial algorithmic method to generate $s$ from $s^*$ (just apply the Turing machine to $s^*$), but
there is no algorithm of length $O(1)$ that computes the shortest description $s^*$ from a general input $s$.
One can show \cite{Vitanyi97} that $s^* \equiv (s,K(s))$. Here, the equivalence symbol  $\equiv$ means that both sides can be obtained from each other
by $O(1)$ programs.
The following equation for the joint algorithmic information of two strings $x,y$  will be  useful \cite{GacsTromp}:
\begin{equation}\label{JointAlg}
K(x,y)\stackrel{+}{=}K(x)+K(y|x^*)=K(x)+K(y|x,K(x))\,.
\end{equation}

The conditional version  reads \cite{GacsTromp}:
\begin{equation}\label{JointAlgCond}
K(x,y|z)\stackrel{+}{=}K(x|z)+K(y|x,K(x|z),z)
\end{equation}

\noindent
The  most important notion in this paper  will be the algorithmic mutual information measuring the amount of
algorithmic information that two objects have in common.
Following \cite{Gruenwald} we  define:

\begin{Definition}[algorithmic mutual  information]${}$\\ 
\label{algomu}
Let $x,y$ be two strings. Then the algorithmic mutual information of   $x,y$ is
\[
I(x:y):=K(y)-K(y|x^*)\,.
\]
\end{Definition}

The mutual information is the number of bits that can be saved in the  description of $y$ when the shortest description of $x$ is
already known.
The fact that one uses $x^*$ instead of $x$ ensures that it  coincides with the
symmetric expression
\cite{Gruenwald}:
\begin{Lemma}[symmetric version of algorithmic mutual information]${}$\\
For two strings $x,y$ we have
\[
I(x:y)\stackrel{+}{=}K(x)+K(y)-K(x,y)\,.
\]
\end{Lemma}

In the following  sections,
non-vanishing mutual  information will be taken as an indicator for causal relations, but
more detailed information on the causal structure will be inferred from {\it  conditional} mutual information.
This is in contrast to approaches from the literature 
to measure similarity versus differences of single objects
that we briefly review now. 
To measure differences between single objects, e.g. pictures
  \cite{BennettDist,LiSim}, one defines the  {\it  information  distance}  $E(x,y)$ between the two corresponding strings
as the  length of  the shortest program that computes $x$ from $y$ {\it and} $y$ from $x$.  
It can be shown \cite{BennettDist} that
\[
E(x,y)\stackrel{{\rm log}}{=} \max\{ K(x|y),K(y|x)\}\,,
\]
where $\stackrel{{\rm log}}{=}$ means equality  up to a logarithmic term. However, whether $E(x,y)$ is small or large is not an appropriate condition for the existence and the  strength of a causal link. Complex objects can have much information in common even though 
 their distance is large. In order  to obtain a measure that relates the amount of information that is disjoint for the two strings
to the amount they  share, Li et al.~\cite{LiSim} and Bennett et al.~\cite{BennettScience} use the ``normalized distance measure''
\[
d_s(x,y):=\frac{K(x|y^*)-K(y|x^*)}{K(x,y)}\Ceq 1-\frac{I(x:y)}{K(x,y)}\,,
\]
or
\[
d(x,y)=\frac{\max\{ K(x|y),K(y|x)\}}{\max\{ K(x),K(y)\}}\,.
\]
The intuitive  meaning of $d_s(x,y)$ is obvious from its direct  relation to mutual information,
and $1-d(x,y)$ measures the fraction of the information of the  more complex string that is shared with the other  one.
Bennett et al. \cite{BennettScience}  propose to construct evolutionary histories of chain letters using such kinds of information distance
measures. However, like in  statistical causal inference, inferring adjacencies on the basis of strongest 
dependences is only possible for simple causal structures like trees.  In the general case, non-adjacent nodes can share more information than
adjacent ones when information is propagated via more than one path. 
Instead  of constructing causal neighborhood relations by comparing
information distances we will therefore use conditional mutual information. 

In order  to define 
its algorithmic  version, we first observe that Definition \ref{algomu}
can be  rewritten into the  less concise form
\[
I(x:y)\stackrel{+}{=}K(y)-K(y|x,K(x))\,.
\]
This formula generalizes more naturally to the conditional analog  \cite{GacsTromp}:

\begin{Definition}[conditional algorithmic mutual information information]${}$\\
\label{Condalgomu}
Let $x,y,z$ be three strings. Then the conditional mutual algorithmic information of   $x,y$, given $z$ is
\[
I(x:y|z):=K(y|z)-K(y|x, K(x|z),z)\,.
\]
\end{Definition}

As shown in \cite{GacsTromp} (Remark II.3), the conditional mutual information also is symmetric up to a constant term:

\begin{Lemma}[symmetric  algorithmic conditional mutual information]${}$\\ \label{symm}
For three strings $x,y,z$ one has:
\[
I(x:y|z)\stackrel{+}{=}K(x|z)+K(y|z)-K(x,y|z)\,.
\]
\end{Lemma}

\begin{Definition}[algorithmic conditional independence]${}$\\ \label{conI}
Given three strings  $x,y,z$, we  call  $x$ conditionally independent of
$y$, given $z$ 
(denoted by  $x\independent y\,|z$) 
if 
\begin{eqnarray*}
I(x:y|z) \approx 0\,.
\end{eqnarray*}
In  words: Given  $z$, the additional 
knowledge of  $y$ does not allow us a stronger compression of $x$. This remains true if we  are given  the 
Kolmogorov complexity of $y$, given $z$.  
\end{Definition}

The theory developed below will describe laws where symbols like $x,y,z$  represent arbitrary strings. Then 
one can always think of {\it sequences} of strings of increasing complexity and statements like ``the equation holds up to constant terms''
are well-defined. We will then understand conditional independence in the sense of  $I(x:y|z)\stackrel{+}{=}0$.
However, if we are talking about three fixed strings that represent objects in  {\it real-life}, this does not make sense 
and
the threshold for considering two strings dependent will heavily depend on the context. For this reason, we will not specify the  symbol $\approx$ any 
further. This is the same arbitrariness as the cutoff rate for statistical dependence tests.

The definitions and lemmas presented so far were strongly motivated by the statistical analog.
Now we want to focus on a theorem in \cite{Gruenwald} that provides a  mathematical relationship
between algorithmic and statistical mutual information. 
First we  rephrase the following theorem 
Theorem 7.3.1 of \cite{Cover}, showing  that
the Kolmogorov complexity of a random string is approximatively given by the entropy of the underlying probability distribution:

\begin{Theorem}[entropy and Kolmogorov complexity]${}$\\ \label{entrKol}
Let ${\bf x}=x_1,x_2,\cdots, x_n$ be a string whose symbols  $x_j \in \cA$ are drawn  i.i.d.  from a probability distribution $P(X)$ 
over the finite alphabet $\cA$. 
Slightly overloading notation, set $P({\bf x}):=P(x_1)\cdots  P(x_n)$. Let $H(.)$ denote the Shannon entropy of a probability distribution. 
Then there is a constant $c$ such that 
\[
H(P(X))\leq \frac{1}{n} E(K({\bf x}|n)) \leq H(P(X))+ \frac{|\cA| \log n}{n} +\frac{c}{n} \quad \forall n\,,
\]
where $E(.)$ is short  hand for the expected value  with respect to $P({\bf x})$.
Hence
\[
\lim_{n\to \infty} \frac{1}{n} E(K(X))=H(P(X))\,.
\]
\end{Theorem}

However,
for our purpose, we need to see the  relation between algorithmic and statistical {\it mutual information}.
If ${\bf x}=x_1,x_2,\cdots, x_n$ and  ${\bf y}=y_1,y_2,\cdots, y_n$   
such that each pair $(x_j,y_j)$  is drawn  i.i.d. from a joint  distribution $P(X,Y)$, the theorem  already shows that
\[
\lim_{n\to \infty} \frac{1}{n} E(I({\bf x}:{\bf y}))=I(X;Y)\,.
\]  
This can be seen by writing statistical mutual information as 
\[
H(P(X))+H(P(Y))-H(P(X,Y))\,.
\]

The above translations between entropy and algorithmic information refer to a particular setting and
to special limits. 
The focus of this paper is mainly the situation where 
the above limits are not justified.  
Before we 
rephrase Theorem 5.3 in \cite{Gruenwald} which provides insights into the general case, we 
recall that a function $f$ is called recursive if there is a program on a Turing machine that 
computes $f(x)$ from the input $x$, and halts on all possible inputs.

\begin{Theorem}[statistical and algorithmic mutual information]${}$
\label{stKmutual}
Given string-valued random variables $X,Y$ with a  recursive probability mass function $P(x,y)$ over pairs $(x,y)$  of strings. 
We then have
\[
I(X;Y)-K(P) \stackrel{+}{\leq} E(I(x:y)) \stackrel{+}{\leq} I(X;Y) +2 K(P)\,,
\]
where $K(P)$ is the length  of the shortest prefix-free program that computes $P(x,y)$ from $(x,y)$.
\end{Theorem}

We want to provide
an intuition about various aspects of this theorem. 

\vspace{0.1cm}
\noindent
(1) If $I(X;Y)$ is large compared  to $K(P)$ the expected algorithmic mutual information is dominated by the
statistical mutual information.  

\vspace{0.1cm}
\noindent
(2) If $K(P)$ is no longer assumed to be small, 
statistical dependences do not necessarily ensure that the  knowledge of $x$ 
allows us to compress $y$ further than without  knowing $x$.  
It could be that the description of the statistical dependences requires more memory
space than its knowledge would save. 

\vspace{0.1cm}
\noindent
(3) 
On the other hand, knowledge  of $x$ could
allow us to compress $y$ even in the case of a product measure on  $x$  and $y$.
Consider,  for instance, the case that we have the point  mass distribution 
on the pair $(x,y)$ with $x=y$.   
To describe a more sophisticated example generalizing this case we first have  to introduce a family of product 
probability distributions on $\{0,1\}^n$ that we will  need several times throughout the paper.

\begin{Definition}[Defining product distributions by strings]${}$\\
\label{StringM}
Let $P_0,P_1$ be two probability distributions on $\{0,1\}$ and $c$ be a binary 
string of length $n$. Then 
\[
{\bf P}_c:=P_{c_1}\otimes  P_{c_2}\otimes \cdots  \otimes P_{c_n}\,
\]
defines a distribution on $\{0,1\}^n$. 
We will later also need  the  following generalization:
If $P_{00},P_{01},P_{10},P_{11}$ are four distributions on $\{0,1\}$, 
then
\[
{\bf P}_{c,d}:=P_{c_1,d_1}\otimes P_{c_2,d_2}\otimes \cdots \otimes P_{c_n,d_n}
\]  
defines also a family of
product measures on  $\{0,1\}^n$ that is labeled by two strings.
\end{Definition}
 
Denote by ${\bf P}_c^{\otimes m}$ the $m$-fold copy of ${\bf P}_c$ from Definition~\ref{StringM}. 
It describes a distribution on $\{0,1\}^{nm}$ assigning
the probbaility ${\bf P}^{\otimes m}_c(x)$ to
  $x\in  \{0,1\}^{nm}$. 
If 
\[
Q(x,y):={\bf P}_c^{\otimes m}(x)\,{\bf  P}_c^{\otimes m}(y)\,,
\]
knowledge of $x$ in the typical case provides knowledge of $c$, provided $m$  is large enough.
Then we can compress  $y$ better than without knowing $x$ because we do not have to describe $c$ any  more.
Hence the algorithmic mutual information  is large and the statistical mutual information is zero because $Q$ is by construction a product distribution. 
In other words, algorithmic dependences  in a setting with i.i.d sampling can arise from statistical dependences
and from algorithmic dependences between probability distributions.

\subsection{Markov  condition for algorithmic dependences among individual objects}

\label{Subsec:Ind}

Now we
state the causal Markov condition for individual objects as a postulate that links algorithmic mutual dependences with causal structure:

\begin{Postulate}[algorithmic causal Markov condition]${}$\\ \label{acMc}
Let $x_1,\dots,x_n$ be $n$ strings representing descriptions of observations 
whose causal connections are formalized by a directed acyclic graph $G$ with $x_1,\dots,x_n$ as nodes. 
Let $pa_j$ be the concatenation of all parents of $x_j$ and  $nd_j$ the concatenation  of all its non-descendants  except $x_j$ itself.
Then 
\[
x_j \independent nd_j \, |pa^*_j\,.
\]  
As in Definition~\ref{conI}, 
the appropriate cut-off rate for  rejecting  $G$ when
$I(x_j:nd_j|pa^*_j)>0$ will not be specified  here.
\end{Postulate}

This formulation is a natural interpretation of Postulate~\ref{abstractMC} in terms of algorithmic independences.
The only point that remains to be justified is why we condition on $pa_j^*$ instead  of $pa_j$,  
i.e., why we are given the  optimal joint compression of
the parent strings. The main reason is  that this turns out to yield nice statements on the equivalence of different
Markov conditions (in analogy to Lemma~\ref{equivMKSt}). 
Since  the differences between $I(x_j:nd_j|pa_j)$ and $I(x_j:nd_j|pa_j^*)$ can only be logarithmic in the string lengths\footnote{this is because $K(x|y)-K(x|y^*)=O(\log |y|)$, see \cite{Vitanyi97}} 
we will  not focus on this issue any further. 

If we apply Postulate~\ref{acMc} to a trivial graph consisting of two unconnected nodes, we obtain 
the following statement.

\begin{Lemma}[causal principle for algorithmic information]${}$\\ \label{CausalPrinciple}
If the  mutual information  $I(x:y)$ between two objects $x,y$  is significantly greater than zero 
they have
some kind of common past. 
\end{Lemma}

Here, common past between two objects means that one has causally influenced the other or there is a third one influencing both.
The statistical version of this principle is part of Reichenbach's principle of the common cause \cite{Reichenbach}
stating that statistical dependences between random variables\footnote{The original formulation  considers actually dependences between events, i.e., binary variables.}
$X$ and $Y$ are always due   at least one of the  following three types of causal  links: (1)  $X$ is a cause of $Y$ or  (2) 
vice versa or (3) there is a common  cause $Z$.  For  objects, the term ``common past'' includes all three types of
causal relations. For a text,  for instance, it reads: similarities of two texts  $x,y$ indicate that one author has been influenced by the other or  that both  have been influenced  by a third one. 

Before we construct a model of  causality that makes it possible to prove the causal Markov  condition we
want to discuss some examples. If one discovers significant similarities in the genome of two sorts of animals one will try to explain  the similarities  by relatedness in the sense of evolution.  Usually, one would, for instance, assume such 
a common history if one has identified
{\it long} substrings that both animals have in common. 
However, the following scenario shows two observations that superficially look similar,  but nevertheless
we cannot infer a common past since their algorithmic complexity is low (implying that the algorithmic mutual information is low, too).

Assume two persons  are instructed to write down a binary string of length $1000$ and both decide to write the  same string 
$x=1100100100001111110...$. It seems straightforward to assume that the persons have communicated
and agreed upon this choice. However, after observing that  $x$ is just the   binary representation of $\pi$, one can easily imagine
that it was just a coincidence that both wrote the same sequence. In other  words, the similarities  are no longer significant after observing
that they stem from a {\it simple} rule. 
This shows that the {\it length} of the pattern that is common to both observations, 
is not a reasonable criterion on whether the similarities are significant.

To understand the  algorithmic causal  Markov condition we will study its implications as well as its  justification. 
In analogy to Lemma~\ref{equivMKSt} we have

\begin{Theorem}[equivalence of algorithmic Markov conditions]${}$ \\ \label{equiAMK}
Given the strings $x_1,\dots,x_n$ and a directed acyclic graph $G$.
Then 
the following conditions are equivalent:

\begin{enumerate}[{\em I.}]

\item {\bf Recursive form:} the joint complexity is given by the
sum of complexities of each node, given the optimal compression of its parents:
\begin{equation}\label{RecAl}
K(x_1,\dots,x_n) \stackrel{+}{=}\sum_{j=1}^n K(x_j|pa^*_j)
\end{equation}

\item {\bf Local Markov condition:} Every node is independent of its non-descendants, given the optimal compression of its parents:
\[
I(x_j:nd_j|pa^*_j) \stackrel{+}{=}0\,.
\]

\item {\bf  Global Markov condition:} 
\[
I(S:T|R^*) \stackrel{+}{=}0
\]
if $R$ d-separates $S$ and $T$.
\end{enumerate}

\end{Theorem}

Below we will therefore no  longer distinguish between  the different versions and just refer to ``the algorithmic Markov condition''.
The intuitive meaning of eq.~(\ref{RecAl}) is that 
the shortest description  of all strings is given by describing how to generate every string from its direct causes.
A similar kind of ``modularity'' of descriptions will also occur later in a different context when we consider description complexity of joint probability 
distributions.

For the proof of Theorem~\ref{equiAMK} we will need a Lemma that 
is an analogue of the observation that for any two random variables $X,Y$ the statistical  mutual information 
satisfies $I(f(X);Y)\leq I(X;Y)$ for every measurable function $f$. 
The algorithmic analog is  to consider two strings $x,y$ and one string $z$ that is derived from $x^*$ by a simple rule.

\begin{Lemma}[monotonicity of algorithmic information]${}$\\ \label{mono}
Let $x,y,z$ be three strings such that $K(z|x^*)\stackrel{+}{=}0$. 
Then
\[
I(z:y)\stackrel{+}{\leq} I(x:y)\,.
\]
\end{Lemma}

This lemma is a special case of  Theorem II.7 in \cite{GacsTromp}.
We will also need the following result:

\begin{Lemma}[monotonicity of conditional information]${}$\\ \label{monC}
Let $x,y,z$ be three strings. Then
\[
K(z|x^*)\stackrel{+}{\geq} K(z|(x,y)^*)\,.
\]
\end{Lemma}

Note that $K(z|x^*)\stackrel{+}{\geq} K(z|x^*,y)$ and $K(z|x^*)\stackrel{+}{\geq} K(z|x^*,y*)$  
is obvious but
 Lemma~\ref{monC} is non-trivial because the star  operation is {\it jointly} applied to $x$ and $y$. 

\vspace{0.3cm}
\noindent
Proof of Lemma~\ref{monC}: Clearly the string $x$ can be derived from $x,y$ by a  program of length  $O(1)$.
Lemma~\ref{mono} therefore implies
\[
I(z:x)\stackrel{+}{\leq} I(z: x,y)\,,
\]
where $I(z:x,y)$  is shorthand for $I(z:(x,y))$.
Hence
\begin{eqnarray*}
K(z)-K(z|x^*)&\stackrel{+}{=}&
 I(z:x) 
\stackrel{+}{\leq} I(z:x,y) \\
&\stackrel{+}{=}& K(z)-K(z|(x,y)^*)\,.
\end{eqnarray*} 
Then we obtain the statement 
by subtracting $K(z)$ and inverting the sign.
$\Box$

\vspace{0.3cm}
The  following lemma will only be used in Subsection~\ref{Subsec:MDL}. We state it here  because it is closely 
related to  the ones above.

\begin{Lemma}[generalized data processing inequality]${}$\\  \label{dataPr}
For any three strings $x,y,z$,
\[
I(x:y|z^*)\stackrel{+}{=}0\,
\]
implies
\[
I(x:y)\stackrel{+}{\leq} I(x:z)\,.
\]
\end{Lemma}

The name ``data processing inequality'' is justified because the assumption $x\independent y\,|z^*$ may arise from
the typical data processing scenario where $y$ is obtained from $x$ via $z$.

\vspace{0.3cm}
\noindent
Proof of Lemma~\ref{dataPr}: 
Using Lemma~\ref{monC} we have
\begin{eqnarray} \label{bcK}
K(x|y^*)&\stackrel{+}{\geq}& K(x|(zy)^*) \\ \nonumber
&\Ceq &K(x|z,y, K(yz))    \\ \nonumber 
&\Ceq  &K(x|z,y ,K(z)+K(y|z^*))\\ \nonumber 
&\stackrel{+}{\geq}& K(x|z,y, K(z),K(y|z^*)\\ \nonumber 
&\Ceq &K(x|z^*,y,K(y|z^*))\,, 
\end{eqnarray}
where the second inequality holds because $K(z)+K(y|z^*)$ can obviously 
be computed from the pair $(K(z),K(y|z^*))$ by an $O(1)$ program. The last  equality
uses, again, the equivalence of  $z^*$ and  $(z,K(z))$.
Hence we obtain:
\begin{eqnarray*}
I(x:y)&\stackrel{+}{=}&K(x)-K(x|y^*) 
\stackrel{+}{=} 
K(x|z^*)+I(x:z)-K(x|y^*)\\
&\stackrel{+}{\leq}&
K(x|z^*)+I(x:z) -K(x|y,K(y|z^*),z^*)\\
 &\stackrel{+}{=}&
I(x:z)+I(x:y|z^*)
\stackrel{+}{=}  I(x:z)\,.
\end{eqnarray*}
The first step is by Definition~\ref{algomu},
the second one uses Lemma~\ref{monC}, the third step is a direct application of ineq.~(\ref{bcK}), the fourth one is due to 
Definition~\ref{Condalgomu}, and the last step is by  assumption. 
$\Box$

\vspace{0.3cm}
\noindent
Proof of Theorem~\ref{equiAMK}: I $\Rightarrow$ III: Define a probability mass function $P$ on 
$(\{0,1\}^*)^{\times n}$, i.e., the set of $n$-tuples of strings, as follows.
Set
\begin{equation}\label{RecKP}
P(x_j|pa_j):=\frac{1}{z_j} 2^{-K(x_j|pa^*_j)}\,,
\end{equation}
where $z_j$ is a normalization factor. In this context, it is important that the symbol $pa_j$ refers to conditioning on
the $k$-tuple of strings $x_i$ that are  parents of $x_j$ (in contrast  to  conditional complexities where we can
interpret $K(.|pa^*_j)$ equally well as conditioning on {\it one} string given by the {\it concatenation} of all those  $x_i$).
 Note that Kraft's inequality (see \cite{Vitanyi97},  Example 3.3.1) implies 
\[
\sum_{x} 2^{-K(x|y)} \leq 1\,,
\]
for every $y$ entailing
that the expression is indeed normalizable by $z_j\leq 1$.
We have
\[
K(x_j|pa^*_j)\stackrel{+}{=}-\log_2 P(x_j|pa_j)\,.
\]
Then we set
\begin{equation}\label{RecStx}
P(x_1,\dots,x_n):=\prod_{j=1}^n P(x_j|pa_j)\,,
\end{equation}
i.e., $P$ is by construction recursive with respect to $G$.
It is easy to see that $K(x_1,\dots,x_n)$ can be computed from $P$:
\begin{eqnarray}\label{RecProof}
K(x_1,\dots,x_n) &\stackrel{+}{=}&\sum_{j=1}^n K(x_j|pa_j^*)  \\ \nonumber
 &\stackrel{+}{=}&-\sum_{j=1}^n \log_2 P(x_j|pa_j)\\   \nonumber
&=&-\log_2 P(x_1,\dots,x_n)\,.
\end{eqnarray}
Remarkably, we can also compute Kolmogorov complexities of {\it subsets} of $\{x_1,\dots,x_n\}$ from the corresponding marginal probabilities.
We start by proving
\begin{equation}\label{cool}
K(x_1,\dots,x_{n-1}) \stackrel{+}{=}-\log_2 \sum_{x_n} 2^{-K(x_1,\dots,x_n)}\,.
\end{equation}
To this end, we observe
\begin{eqnarray}\label{upperMarg}
\sum_{x_n} 2^{-K(x_1,\dots,x_n)} &\stackrel{\times}{=}& \sum_{x_n} 2^{-K(x_1,\dots,x_{n-1})-K(x_n|(x_1,\dots,x_ {n-1})^*)} \\ 
&\stackrel{\times}{\leq}&
2^{-K(x_1,\dots,x_{n-1})}\,, \nonumber
\end{eqnarray}
where $\stackrel{\times}{=}$ denotes equality up to a multiplicative constant. The equality follows from eq.~(\ref{JointAlg}) and
the inequality is obtained by applying Kraft's inequality \cite{Vitanyi97}  to the conditional complexity $K(.|(x_1,\dots,x_ {n-1})^*)$.
On the other hand we have 
\[
K(x_1,\dots, x_{n-1}) \stackrel{+}{=}K(x_1,\dots,x_{n-1},0)\,,
\]
since adding the $1$-bit string $x_n=0$ certainly can be performed by a program of length $O(1)$.
Hence we have
\begin{eqnarray*}
\sum_{x_n} 2^{-K(x_1,\dots,x_n)} &\stackrel{\times}{\geq}&2^{-K(x_1,\dots,x_{n-1},0)}  \\&\stackrel{\times}{=}& 2^{-K(x_1,\dots,x_{n-1})}  \,.
\end{eqnarray*}
Combining this with ineq.~(\ref{upperMarg}) yields
\[
2^{-K(x_1,\dots,x_{n-1})} \stackrel{\times}{=} \sum_{x_n} 2^{-K(x_1,\dots,x_n)}\,.
\]
Using eq.~(\ref{RecProof}) we obtain
\begin{eqnarray*}
K(x_1,\dots,x_{n-1})&\stackrel{+}{=}&-\log_2 \sum_{x_n} 2^{-K(x_1,\dots,x_n)}\\ &\stackrel{+}{=}&-\log_2 \sum_{x_n} P(x_1,\dots,x_n)\\ &\stackrel{+}{=}&-\log_2 P(x_1,\dots,x_ {n-1})\,,
\end{eqnarray*}
which proves equation~(\ref{cool}).
This implies
\[
K(x_1,\dots,x_ {n-1}) \stackrel{+}{=}-\log_2  P(x_1,\dots,x_{n-1})\,.
\]
Since the same argument holds for marginalizing over any other variable $x_j$ we conclude that
\begin{equation}\label{subK}
K(x_{j_1},\dots,x_{j_k}) \stackrel{+}{=} -\log_2 P(x_{j_1},\dots,x_{j_k})\,,
\end{equation}
for every subset of strings of size $k$ with $k\leq n$. This follows by induction over $n-k$.

Now we can
use  the relation between marginal probabilities and Kolmogorov complexities to show that conditional complexities are also given by the corresponding 
{\it conditional} probabilities, i.e., for any two subsets $S,T\subset \{x_1,\dots,x_n\}$ we have
\[
K(S|T^*)\stackrel{+}{=}-\log_2 P(S|T)\,.
\]
Without loss of generality, set $S:=\{x_1,\dots,x_j\}$ and $T:=\{x_{j+1},\dots,x_k\}$ for $j<k\leq n$.
Using  eqs.~(\ref{JointAlg}) and (\ref{subK}) we get 
\begin{eqnarray*}
K(x_1,\dots,x_j|(x_{j+1},\dots,x_k)^*)&\stackrel{+}{=}&K(x_1,\dots,x_k)-K(x_{j+1},\dots,x_k)\\ &\stackrel{+}{=}&-\log_2 \Big( P(x_1,\dots,x_k)/P(x_{j+1},\dots,x_k) \Big)
\\&\stackrel{+}{=}&
-\log_2 P(x_1,\dots,x_j|x_{j+1},\dots,x_k)\,.
\end{eqnarray*}
Let $S,T,R$ be three subsets of $\{x_1,\dots,x_n\}$ 
such that $R$ d-separates $S$ and $T$. 
Then
$S\independent T\,|R$ with respect to $P$ because $P$ satisfies the recursion~(\ref{RecStx})
(see Lemma~\ref{equivMKSt})\footnote{Since $P$ is, by construction, a discrete probability function, $P$ the density with respect to a  product measure is  directly given by  the probability function itself.}.
Hence 
\begin{eqnarray*}
K(S,T|R^*)&\stackrel{+}{=}&-\log_2 P(S,T|R) \\&\stackrel{+}{=}&
-\log P(S|R)-\log_2 P(T|R)\\ &\stackrel{+}{=}& 
K(S|R^*)+K(S|R^*)\,.
\end{eqnarray*}
This proves algorithmic independence of $S$ and  $T$, given $R^*$ and thus  I $\Rightarrow$ III.
 
To show that III $\Rightarrow$ II it suffices to recall  that $nd_j$ and $x_j$ are d-separated by $pa_j$.
Now we show II $\Rightarrow$ I in strong analogy to the proof for the statistical version of this statement in \cite{LauritzenDawid}:
Consider first a terminal node of $G$. Assume, without loss of generality, that it is $x_n$. 
Hence all strings $x_1,\dots,x_{n-1}$ are non-descendants of $x_n$. 
We thus have $(nd_n,pa_n) \equiv (x_1,\dots,x_{n-1})$ where $\equiv$ means that both strings coincide up to a permutation (on one side) and removing
those strings that occur twice (on the other side).
Due to eq.~(\ref{JointAlg}) we have
\begin{equation}\label{RecAn}
K(x_1,\dots,x_n) \stackrel{+}{=}K(x_1,\dots,x_{n-1})+K(x_n|(nd_n,pa_n)^*)\,.
\end{equation}
Using, again,  the equivalence of $w^*\equiv (w,K(w))$ for any string  $w$   we have
\begin{eqnarray}
K(x_n|(nd_n,pa_n)^*)&\stackrel{+}{=}& K(x_n|nd_n,pa_n, K(nd_n,pa_n)) \nonumber \\ \nonumber
&\stackrel{+}{=}& K(x_n|nd_n,pa_n,K(pa_n)+K(nd_n|pa^*_n))\\  \nonumber &\stackrel{+}{\geq} & 
K(x_n|nd_n,pa^*_n,K(nd_n|pa^*_n))\\  &\stackrel{+}{=}& 
K(x_n|pa_n^*)\,.  \label{vpan} 
\end{eqnarray}
The second step follows from $K(nd_n,pa_n)=K(pa_n)+K(nd_n|pa_n^*)$. The inequality holds because 
$nd_n,pa_n,K(pa_n)+K(nd_n|pa_n^*)$ can be computed from $nd_n,pa_n^*,K(nd_n|pa^*_n)$ via a program of length  $O(1)$.
The last step follows directly from the assumption  $x_n \independent nd_n \,| pa_n^*$. 
Combining ineq.~(\ref{vpan}) with  Lemma~\ref{monC} yields
\begin{equation}\label{upperxn}
K(x_n|(nd_n,pa_n)^*)\stackrel{+}{=} K(x_n|pa^*_n)\,.
\end{equation}
Combining eqs.~(\ref{upperxn}) and  (\ref{RecAn}) we obtain
\begin{equation}
K(x_1,\dots,x_n) \stackrel{+}{=}K(x_1,\dots,x_{n-1})+K(x_n|pa_n^*)\,.
\end{equation}
Then statement I follows by induction over $n$. 
$\Box$

\vspace{0.3cm}
\noindent
To show  that the algorithmic Markov condition can be derived from an algorithmic version of 
the functional model in  Postulate~\ref{functional} we 
introduce  the following model of  causal mechanisms.

\begin{Postulate}[algorithmic model of causality]${}$\\  \label{algoFunc}
Let $G$ be a DAG formalizing the causal structure among the strings $x_1,\dots,x_n$.
Then every $x_j$ 
is computed by a program $q_j$ with length  $O(1)$ from its parents $pa_j$ and an additional input $n_j$. We write formally
\[
x_j=q_j(pa_j,n_j)\,,
\]
meaning that the Turing machine computes $x_j$ from the input $pa_j,n_j$ using the additional program  $q_j$ and halts.
The inputs $n_j$ are jointly independent in the sense
\[
n_j  \independent n_1,\dots,n_{j-1},n_{j+1},n_n\,.
\] 
By defining new programs  that contain $n_j$ we  can, equivalently, drop the assumption that
the programs $q_j$ are simple and assume  that  they are  jointly independent instead.
\end{Postulate}

We could also have assumed that $x_j$ is a function  $f_j$ of all its parents, but our model is more general
since the map defined by the input-output behavior of $q_j$ need not be a total function \cite{Vitanyi97}, i.e.,  the Turing machine 
simulating the process would not necessarily halt
on {\it all} inputs $pa_j,n_j$. 

The  idea  to represent causal mechanisms by programs written for some universal 
Turing machine is basically in the spirit of various interpretations of the Church-Turing thesis.
One formulation, given by Deutsch \cite{Deutsch85}, states that every process 
taking place
in the real world can  be simulated  by a Turing machine. Here we  assume that 
the way different systems influence each other by  physical signals can be  simulated by computation  processes that exchange 
messages  of bit strings.\footnote{Note, however, that 
sending quantum systems between the nodes could 
transmit a kind of information (``quantum information'' \cite{NC}) that cannot  be phrased  in terms  of bits. 
It is known that this enables completely new communication scenarios,  e.g. quantum cryptography.  
The relevance 
of quantum information transfer for  causal inference is not yet fully understood. It has, for instance,
been shown that the violation of Bell's inequality in quantum  theory is also relevant for causal inference \cite{Drugs}.
This is because some causal inference rules   
between classical variables break down when the latent factors are represented by {\it quantum} states rather than being 
classical variables.}

Note that mathematics also allows us  to construct strings that are linked to each other in an {\it uncomputable} way. 
For instance, let $x$ be an arbitrary binary string 
and  $y$ be defined by $y:=K(x)$.  However, it is hard to believe that a real causal mechanism could create such kind of 
relations between
objects given that one believes that real processes can always be simulated by algorithms. These remarks are intended
to  give sufficient motivation for our model.

Postulate~\ref{algoFunc} implies the algorithmic causal Markov condition:

\begin{Theorem}[algorithmic model implies Markov]${}$\\ \label{algoImplM}
Let $x_1,\dots,x_n$ be generated by the model in Postulate \ref{algoFunc}.
Then they satisfy the algorithmic Markov condition with respect to $G$. 
\end{Theorem}

\noindent
Proof (straightforward adaption of the proof of Lemma~\ref{funcM}):
Extend $G$ to a causal structure $\tilde{G}$ with nodes 
$x_1,\dots,x_n,n_1,\dots,n_n$.   
To see that the extended set of nodes satisfy the local Markov condition w.r.t. $\tilde{G}$, observe first  that 
every node $x_j$ is given by its parents via an $O(1)$ program. Second, every $n_j$ is parentless and (unconditionally)
independent of all its non-descendants because they can be  computed from $\{n_1,\dots.n_n\}\setminus \{n_j\}$ 
via an  $O(1)$ program.

By Theorem~\ref{equiAMK} the extended set  of  nodes is also globally Markovian w.r.t. $\tilde{G}$. 
The parents $pa_j$ d-separate
$x_j$ and $nd_j$ in $\tilde{G}$ (here the parents $pa_j$ are still defined  with respect to $G$). 
This implies the local Markov condition for $G$.  
$\Box$

\vspace{0.3cm}
\noindent
It is trivial to  construct examples where the causal  Markov condition is violated if the programs
$q_j$ are mutually dependent (for instance, the trivial graph with  two nodes $x_1,x_2$ and no edge would satisfy
$I(x_1:x_2)>0$ if the  programs $q_1,q_2$  computing $x_1,x_2$ from an empty input are dependent).

The last sentence of Postulate~\ref{algoFunc} makes apparent that 
the {\it mechanisms} that  generate causal relations are assumed to be independent.
This 
is essential for the general philosophy  of this paper.
To see that such a mutual independence of mechanisms is a reasonable  assumption we recall that the causal graph
is meant to formalize {\it all} relevant causal links between the objects. If we observe, for instance, that two nodes are generated 
from their parents by the same complex rule we postulate another causal link between the nodes that explains the similarity of
mechanisms.\footnote{One could argue that this would be just the causal principle implying 
that similarities of the ``machines'' generating $x_j$ from $pa_j$ 
has to be explained by a causal relation, i.e., a common past of the machines. However, in the context of this paper, such
an argument would be  circular. We have argued that the causal principle is a special case of  the Markov condition and 
derived the latter from the algorithmic model above. We will therefore consider the independence of mechanisms as a first principle.}

\subsection{Relative causality}

\label{Subsec:Rel}

This subsection explains why
it is sensible to define algorithmic dependence and the existence or non-existence of  causal links
{\it relative} to some background information. 
To  this end, we consider genetic sequences  $s_1,s_2$ of two persons that are not relatives. 
We certainly find high 
similarity that leads to  a significant violation of $I(s_1:s_2)= 0$ due  to the fact that both genes are taken from humans. 
However, given the background information ``$s_1$ is a human  genetic sequence'', $s_1$ can be further compressed. 
The same applies to $s_2$.   
Let $h$ be a  code that is particularly adapted to the human genome  in the sense that 
the expected conditional Kolmogorov complexity, given $h$, of a randomly chosen human genome is minimal. 
Then it would make sense to consider $I(s_1:s_2|h)>0$ as a hint for a relation that goes beyond the fact that both persons are human.
In contrast, for  the unconditional mutual information we expect $I(s_1:s_2) \geq K(h)$. 
We will 
therefore infer some causal relation (here: common ancestors in the evolution) using the causal principle in Lemma~\ref{CausalPrinciple}
(cf.~\cite{Milo}).

The  common  properties between different and unrelated individuals  of the same species can be screened off by providing the relevant background information. Given this causal  background, we can detect further similarities in the genes by the conditional  algorithmic mutual information and take them as an indicator for an additional causal relation that goes beyond the common evolutionary background. 
For this reason, every discussion on whether there exists a causal link between two objects (or individuals) requires
a specification of the background information. In this sense, causality is a relative concept. 

One may ask whether such a relativity of causality is also true for the statistical version of the causality principle, i.e.,
Reichenbach's principle of the common cause. 
In the statistical version of the link between causality and dependence, the  relevance of the background information is 
less obvious because it is evident that statistical  methods are always applied to a {\it given statistical  ensemble}.
If we, for instance, ask whether  there is a  causal relation between the height and the income of a person without specifying 
whether we   refer to people of a certain age, we observe the same relativity with respect to additionally
specifying the ``background information'', which is  here given by referring to a specific ensemble.  

In the following sections we will assume that the relevant background information has been specified and it has been clarified how to translate the  relevant aspects of a real object into a binary  string such that we  can identify every object with its  binary description.

\section{Novel statistical inference rules from the algorithmic Markov condition}

\label{Sec:NewInf}

\subsection{Algorithmic independence of Markov kernels}

\label{Subsec:plMKAlg}

To describe the implications of  the algorithmic Markov  condition for statistical causal inference, 
we consider random variables $X$ and $Y$ where $X$ causally influences $Y$.  We can think  of  $P(X)$ as describing a source  $S$
that generates $x$-values and sends them to a ``machine''   $M$ that generates $y$-values according to $P(Y|X)$.
Assume we observe that
\[
I(P(X):P(Y|X))\gg 0\,.
\]
Then we conclude that there must be a causal link between $S$ and $M$ that  goes beyond transferring $x$-values from $S$ to $M$. This is because $P(X)$ and $P(Y|X)$ are inherent properties of $S$ and $M$, respectively which do not  depend on the current value
of $x$ that has been sent.
Hence there must be a  causal link that explains the similarities in the {\it design} of  $S$ and $M$.
Here we have assumed that we know that $X\rightarrow Y$ is the correct causal  structure on the statistical level. Then we have  to accept
that a causal link on the level of the machine design is present.

If the causal structure on the statistical level is unknown, we would prefer causal hypotheses that explain the data without needing a
causal connection  on the higher level provided that they satisfy the statistical  Markov condition.
Given this principle, we thus will  prefer causal graphs $G$ for which the Markov kernels
$
P(X_j|PA_j)$ become algorithmically independent. 
This is equivalent to saying that the shortest description of $P(X_1,\dots,X_n)$ is given by concatenating the descriptions of the Markov
kernels, a postulate that has already been formulated by Lemeire and Dirkx \cite{LemeireD}:

\begin{Postulate}[algorithmic independence of statistical properties]${}$\\ \label{inMech}
A causal hypothesis $G$ (i.e., a DAG) is only acceptable if
the shortest description of the joint density $P$ is given by a concatenation
of the shortest description of the Markov kernels, i.e.
\begin{equation}\label{Ksum}
K(P(X_1,\dots,X_n))\Ceq \sum_j K(P(X_j|PA_j))\,.
\end{equation}
If no such causal graph exists, we reject  every possible DAG and assume that
there is a causal relation of a different type, e.g.,
a latent common cause, selection bias,  or
a cyclic causal structure.  
\end{Postulate}

The sum on the right hand side of eq.~(\ref{Ksum}) will be called the {\it total complexity} of the causal model $G$.
Note that Postulate~\ref{inMech} implies that we have to reject every causal hypothesis for which
the total complexity is not minimal because  a model with shorter total complexity already provides a shorter description of the
joint distribution. 
Inferring causal directions by
minimizing this expression (or actually a computable modification) could also be  interpreted in a Bayesian way if we consider $K(P(X_j|PA_j))$ 
as the negative log likelihood for the prior probability for having the conditional $P(X_j|PA_j)$ (after appropriate normalization). 
However, Postulate~\ref{inMech} contains an idea that goes beyond known Bayesian approaches to causal discovery because it provides hints on the incompleteness
of the class of models under consideration (in addition to providing rules for giving preference {\it within} the class).

Lemeire and Dirkx
\cite{LemeireD} already show that the causal faithfulness principle (Postulate~\ref{faithfulness}) follows
from Postulate~\ref{inMech}.
Now we want  to show that it also implies causal inference rules  that go beyond
the known ones. 

To this end, we  focus again on the example in Subsection~\ref{Subsec:plMK} with a binary variable  $X$ and a continuous variable $Y$.
The  hypothesis $X\rightarrow Y$ is not  rejected on the basis of Postulate~\ref{inMech} 
because  $I(P(X):P(Y|X))\Ceq 0$. For  the equally weighted mixture of  two Gaussians this already  follows\footnote{for the more general case $P(X=1)=p$ with $K(p)\gg 0$, 
this also follows if we assume that $p$ is algorithmically independent of the
parameters that specify  $P(Y|X)$}  from $K(P(X))\Ceq 0$.
On the other hand, $Y\rightarrow X$ violates Postulate~\ref{inMech}.
Elementary calculations show that the conditional $P(X|Y)$ is  given by
the sigmoid function
\begin{equation}\label{sig}
P(X=1|y)=\frac{1}{2}\Big( 1+\tanh \frac{\lambda (y-\mu)}{\sigma^2} \Big)\,.
\end{equation}
We observe  that the
same parameters $\sigma,\lambda,\mu$ that occur  in  $P(Y)$, also occur in $P(X|Y)$. This already shows that
the two Markov kernels are algorithmically dependent. 
To be more explicit, we observe that $\mu$, $\lambda$, and  $\sigma$ are required
 to specify  $P(Y)$. To describe $P(X|Y)$, we need $\lambda/\sigma^2$ and $\mu$. 
Hence we have
\begin{eqnarray*}
K(P(Y))\Ceq K(\mu,\lambda,\sigma) &\Ceq& K(\mu)+K(\lambda)+K(\sigma)\\
K(P(X|Y))\Ceq K(\mu,\lambda/\sigma^2) &\Ceq & K(\mu)+K(\lambda/\sigma^2)\\
K(P(X,Y))\Ceq K(P(Y),P(X|Y))&\Ceq& K(\mu,\lambda,\sigma)\Ceq K(\mu)+K(\lambda)+K(\sigma)\,,
\end{eqnarray*}
where we have assumed that the strings $\mu,\lambda,\sigma$ are jointly independent. 
Note that the information that $P(Y)$ is a mixture of two  Gaussians
and that $P(X|Y)$ is a  sigmoid counts as a constant because its description complexity 
does not depend on  the parameters.  

We  thus get
\[
I(P(Y):P(X|Y))\Ceq K(\mu)+K(\lambda/\sigma^2)\,.
\]
Therefore we reject the causal  hypothesis  $Y\rightarrow X$ due to Postulate~\ref{inMech}.
The  interesting point  is that 
we need not look at the alternative hypothesis $X\rightarrow Y$. In other words, 
we do not reject $Y\rightarrow X$ {\it only} because the converse direction leads to simpler expressions. We can reject 
it alone one the basis of observing algorithmic dependences between $P(Y)$ and  $P(X|Y)$ making the causal model suspicious. 

\begin{figure}
\centerline{\includegraphics[scale=0.25]{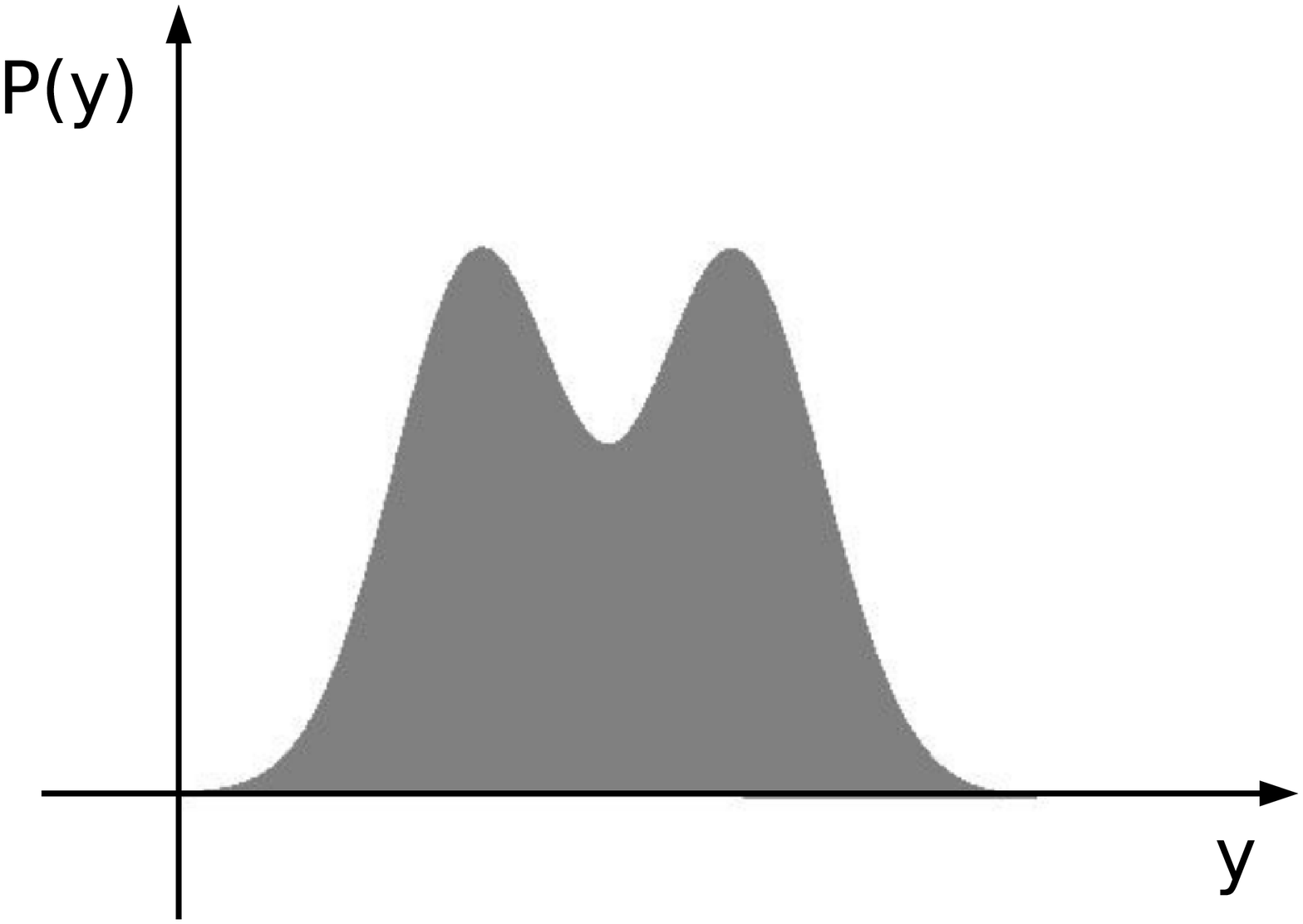}
\includegraphics[scale=0.25]{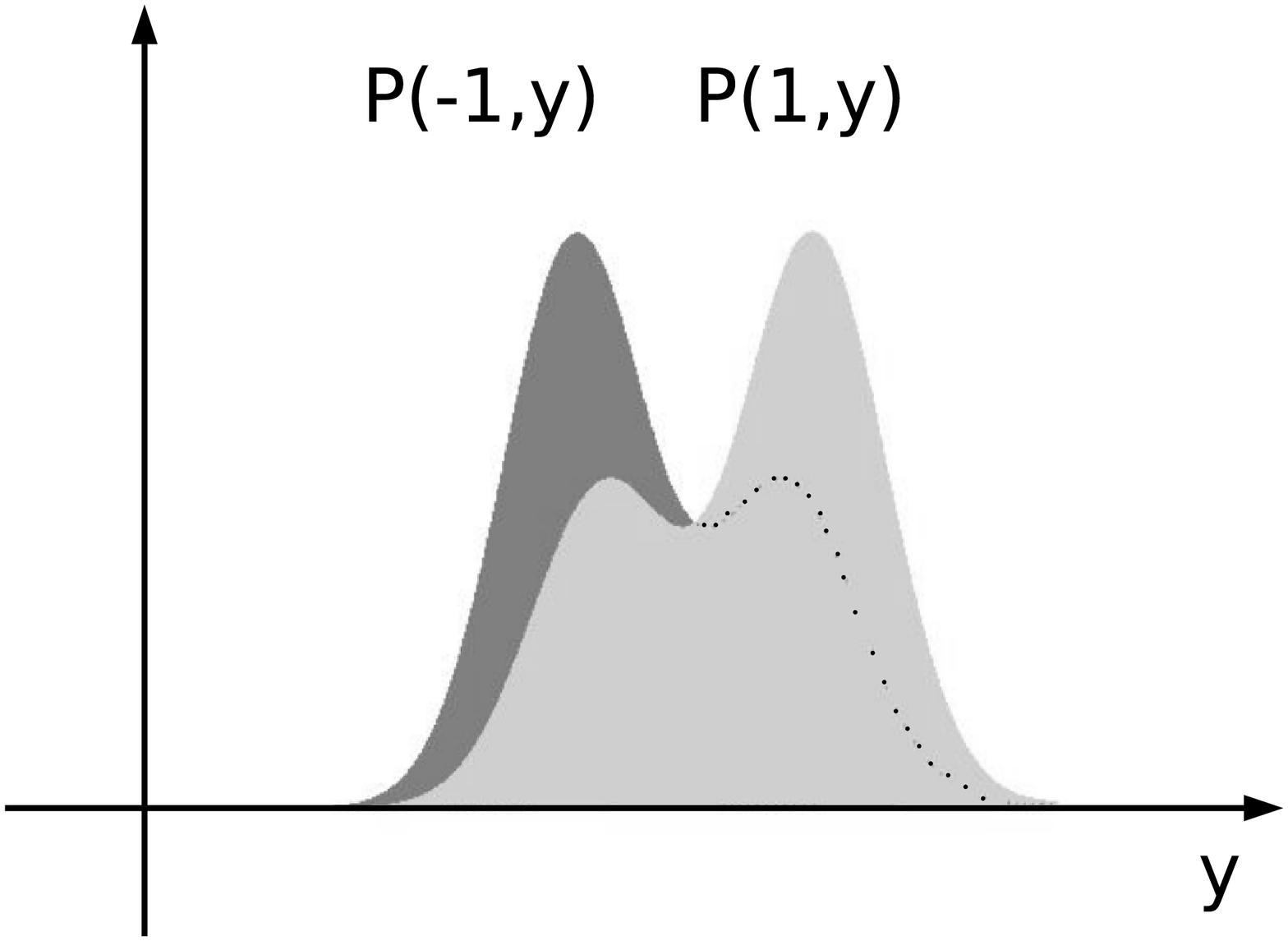}}
\caption{{\small Left: a source generates the bimodal distribution $P(Y)$. A  machine generates $x$-values according to 
a conditional $P(X|Y)$  given by a sigmoid function. If  the slope and the position parameters of the sigmoid are not correctly 
adjusted to the 
distance, the position, and the width of the two Gaussian modes, 
the generated joint distribution no longer consists of two Gaussians (right).}\label{GaussianDet}} 
\end{figure}

The following gedankenexperiment shows that $Y\rightarrow X$ would become plausible if we 
``detune'' the sigmoid $P(X|Y)$ by choosing 
$\tilde{\lambda},\tilde{\mu},\tilde{\sigma}$ independently of $\lambda$ and $\mu$, and $\sigma$. 
Then $P(Y)$ and $P(X|Y)$ are by definition  algorithmically independent and therefore we obtain a more complex joint distribution:
\[
K(P(X,Y))=K(\lambda)+K(\tilde{\lambda})+K(\mu)+K(\tilde{\mu})+K(\sigma)+K(\tilde{\sigma})\,. 
\]  
The fact that the set of mixtures  of two Gaussians does not have six free parameters already shows that 
$P(X,Y)$ must be a more complex distribution than the one above.   
Fig.~\ref{GaussianDet} shows an example of a joint distribution obtained for the  ``detuned''  situation.

As  already noted by \cite{LemeireD},  the independence of mechanisms is related to Pearl's thoughts on the stability of causal statements: 
the causal mechanism $P(X_j|PA_j)$ does not change  if one changes  the input distribution $P(PA_j)$
by influencing the variables $PA_j$. The same conditional  can therefore
occur, under different background conditions, with different input distributions. 

Postulate~\ref{inMech} naturally occurs in  the probability-free
version of the causal Markov condition. 
To explain this, assume we are given two strings ${\bf x}$ and ${\bf y}$ of length  $n$
(describing two real-world observations)
and noticed that  ${\bf x}={\bf y}$.  Now we consider two alternative scenarios:

\vspace{0.2cm}
\noindent
(I) Assume that every pair $(x_j,y_j)$ of digits ($j=1,\dots,n$) has 
been independently drawn from the same joint distribution $P(X,Y)$ of the binary random variables $X$ and $Y$.

\vspace{0.2cm}
\noindent
(II)
Let ${\bf x}$ and ${\bf y}$ be single instances of
string-valued random variables $X$ and $Y$. 

\vspace{0.2cm}
The  difference between (I)  and (II) is crucial for statistical causal inference:
In case (I), statistical  independence is rejected with high confidence proving the existence of a causal link.
In constrast,
there is  no evidence for statistical dependence in case (II) since the underlying 
joint distribution
on $\{0,1\}^n\times \{0,1\}^n$ could, for instance, be the point mass on the pair  $({\bf x},{\bf y})$, 
which is a product distribution, i.e.,
\[
P(X,Y)=P(Y)P(X)\,.
\]
Hence,
statistical causal inference would not infer a causal connection in case (II).

Algorithmic causal inference, on the other hand, infers a causal link in both cases because
the equality ${\bf x}={\bf y}$ requires an explanation.
The relevance of switching  between (I) and (II) then consists merely in shifting the 
causal connection
to another level: In the i.i.d setting,  every $x_j$ must be causally linked to $y_j$.
In case (II), 
there must be a connection between the two {\it mechanisms} that have generated
the entire strings because $I(P(X):P(Y|X))=I(P(X):P(Y))\gg 0$.
This can, for instance, be  due to the fact that two machines emitting the same string were designed
by the same engineer. A detailed discussion of the relevance of translating the i.i.d.  assumption into the setting of 
algorithmic causal inference will be given in
Subsection~\ref{Subsec:Res}.

\subsection*{Examples with large probability spaces}


In the preceding subsection we have ignored a serious problem  with defining the Kolmogorov complexity of (conditional) probability
distributions that even occurs in finite  probability spaces.  
First of all the ``true'' probabilities may not be computable. For instance, a coin may produce ``head'' with  
probability $p$ where $p$ is some {\it uncomputable} number, i.e., $K(p)=\infty$.
And even if  it were some  computable value $p$ with large $K(p)$ it would  be quite artificial to call the probability
distribution $(p,1-p)$ ``complex'' because $K(p)$ is high and ``simple'' if we have,  for instance $p=1/\pi$. 
 A more reasonable notion of complexity 
can be obtained by describing the probabilities only up to a certain accuracy $\epsilon$. If $\epsilon$ is not to small
we obtain small complexity values for the  distribution  of a binary variable, and also low complexity for a
distribution on a larger set that is $\epsilon$-close to the values of some simple analytical expression like a Gaussian distribution.
There will still remain some unease about the concept of Kolmogorov complexity of ``the true distribution''. We will subsequently develop a  
formalism that avoids this concept. However, Kolmogorov complexity of distributions is a useful idea to start with since it provides
an intuitive understanding of the roots of the asymmetries between cause and  effects 
that we will describe  in Subsection~\ref{Subsec:Res}.

Below, we  will describe a gedankenexperiment with two random variables $X,Y$ linked by the causal structure $X\rightarrow Y$  
where the total complexities of  the  
causal models $X\rightarrow  Y$ and $Y\rightarrow X$ both are well-defined and, in the generic case, different.  
First we will show that they can at most differ by factor two.

\begin{Lemma}[maximal complexity quotient]${}$\\ \label{2bound}
For every joint  distribution $P(X,Y)$  we have
\[
K(P(Y))+K(P(X|Y))\stackrel{+}{\leq} 2 \Big(  K(P(X))+ K(P(Y|X))\Big)\,.
\]
\end{Lemma}

\noindent
Proof: 
Since marginals and conditionals both can be computed from $P(X,Y)$ we have
\[
K(P(Y))+K(P(X|Y))\stackrel{+}{\leq} 2 K(P(X,Y))\,.
\]
Then the statement follows  because $P(X,Y)$ can be computed from $P(X)$ and  $P(Y|X)$.
$\square$

\vspace{0.3cm}
\noindent
To construct  examples  where the bound in Lemma~\ref{2bound} is attained we first introduce a 
method to construct conditionals with well-defined complexity:

\begin{Definition}[Conditionals and joint distributions from strings]\label{StringJ}${}$\\
Let $M_0,M_1$ be two stochastic matrices that  specify transition probabilities  from $\{0,1\}$ to $\{0,1\}$.
Then
\[
{\bf M}_c:=M_{c_1}\otimes M_{c_2}\otimes \cdots \otimes M_{c_n}\,
\]
defines transition probabilities  from $\{0,1\}^n$ to $\{0,1\}^n$. 

We also introduce the same construction for double indices:
Let $M_{00},M_{01},M_{10},M_{11}$ be stochastic matrices describing transition probabilities from  $\{0,1\}$ to $\{0,1\}$.
Let $c,d\in \{0,1\}^n$ be two strings. Then
\[
{\bf M}_{c,d}:=M_{c_1,d_1}\otimes M_{c_2,d_2}\otimes \cdots \otimes M_{c_n,d_n}
\]
defines a transition matrix  from $\{0,1\}^n$ to $\{0,1\}^n$. 
If the matrices $M_j$ or $M_{ij}$ denote joint distributions on  $\{0,1\}\times\{0,1\}$
the objects ${\bf  M}_{c}$ and ${\bf M}_{c,d}$ define joint distributions on $\{0,1\}^n\times \{0,1\}^n$ 
in a canonical way.
\end{Definition}

Let $X,Y$ be variables whose values are the set  of strings in $\{0,1\}^n$. 
Define distributions $P_0,P_1$ on $\{0,1\}$ and stochastic matrices $A_0,A_1$ describing transition  probabilities from  
$\{0,1\}$ to $\{0,1\}$.
Then a string $c\in \{0,1\}^n$ defines a distribution  
$P(X):={\bf  P}_c$ (using Definition~\ref{StringM}) that has well-defined Kolmogorov complexity $K(c)$ if the description 
complexity of $P_0$ and $P_1$ is neglected. 
Furthermore, we set $P(Y|X):={\bf A}_d$ as in Definition~\ref{StringJ}, 
where we have used the canonical identification between stochastic matrices and conditional probabilities
and $d\in \{0,1\}^n$ denotes  some randomly chosen string. 
Let $R_{ij}$  denote    
the joint distribution on $\{0,1\}\times \{0,1\}$ induced  by the marginal $P_i$ on the first  component and the conditional
$A_j$ for  the second, given the first.
Denote the corresponding marginal  dsitribution on the right component by $Q_{ij}$, i.e.,
\[
Q_{ij}:=A_j P_i\,,
\]
and let $B_{ij}$ be the stochastic matrix that describes the conditional probability for the first component, given the second.

Using these notations and the ones in Definition~\ref{StringJ}, we  obtain
\begin{eqnarray}
P(X)&=&{\bf P}_c  \label{tableCM}\\ \nonumber
P(Y|X)&=&{\bf A}_{d}\\ \nonumber 
P(X,Y)&=&{\bf R}_{c,d}\\ \nonumber 
P(Y)&=& {\bf Q}_{c,d}\\ \nonumber
P(X|Y)&=&{\bf B}_{c,d} 
\end{eqnarray}

It is noteworthy  that $P(Y)$ and  $P(X|Y)$ are labeled by both strings while $P(X)$ and  $P(Y|X)$ are described  by only one string  each.
This already suggests that the latter are more complex in the generic case.

Now we compare the sum $K(P(X))+K(P(Y|X))$ to $K(P(Y))+K(P(X|Y))$ for the case$K(c)\Ceq K(d) \Ceq n$.
We assume that  $P_i$ and $A_j$ are computable and  their complexity is counted as $O(1)$ because it does  not depend on $n$. 
 Nevertheless, we assume that $P_i$ and $A_j$ are ``generic'' in the following sense:
All marginals $Q_{ij}$ and conditionals $B_{ij}$ are different whenever $P_0\neq P_1$ and $A_0\neq A_1$.
If we impose one  of the conditions $P_0=P_1$ and $A_0= A_1$   or both, we assume that only those
marginals $Q_{ij}$ and conditionals $B_{ij}$ coincide for which the equality follows from the conditions imposed.
Consider the  following cases:
\vspace{0.3cm}

\noindent
{\bf Case 1:} $P_0=P_1$, $A_0=A_1$.
Then all the complexities vanish because the joint distribution does not depend on the strings $c$ and $d$.

\vspace{0.2cm}
\noindent
{\bf Case 2:} $P_0\neq P_1$, $A_0=A_1$. Then the digits of $c$  are relevant, but the digits of $d$ are not.
Those marginals and conditionals in table~(\ref{tableCM}) that formally depend
on $c$ and $d$, as well as those that depend on $c$, have complexity $n$. Those depending on $d$ have complexity $0$.
\[
\begin{array}{cccccc}
K(P(X))+K(P(Y|X))&\Ceq &n+0&=&n& \\
K(P(Y))+K(P(X|Y))&\Ceq &n+n&=&2n&.
\end{array}
\]

\vspace{0.2cm}
\noindent
{\bf Case 3:} $P_0 =P_1$, $A_0\neq A_1$. Only the dependence on $d$ contributes to the complexity.
This implies
\[
\begin{array}{cccccc}
K(P(X))+K(P(Y|X))&\Ceq&0+n&=&n&\\
K(P(Y))+K(P(X|Y))&\Ceq &n+n&=&2n&.
\end{array}
\]

\vspace{0.2cm}
\noindent
{\bf Case 4:} $P_0\neq P_1$ and $A_0\neq A_1$. Every formal dependence  of the conditionals and marginals on $c$ and $d$ in table~(\ref{tableCM})
is a proper dependence.
Hence we obtain
\[
\begin{array}{cccccc}
K(P(X))+K(P(Y|Y))&\Ceq  &n+n&=&2n&\\
K(P(Y))+K(P(X|Y))&\Ceq &2n+2n&=&4n&.
\end{array}
\]

\vspace{0.3cm}
The general principle of the above example is very  simple. Given that $P(X)$ is taken from  a model class that consists of $N$ different elements and
$P(Y|X)$ is  taken from a  class with $M$ different elements. Then the class of possible $P(Y)$ and the class of possible $P(X|Y)$ both can contain $N\cdot M$ elements.
If  the simplicity of a model is quantified in terms of the size of the class it is taken from (within a hierarchy of more and more complex models),
the statement that $P(Y)$ and $P(X|Y)$ are typically complex is just based on this simple counting argument.

\subsection*{Detecting common causes via dependent Markov kernels}

The following  model shows that latent common causes can yield joint distributions whose Kolmogorov  complexity is smaller 
than $K(P(X))+K(P(Y|X))$ and $K(P(Y))+K(X|Y))$.
Let $X,Y,Z$ have values in $\{0,1\}^n$ and let $P(Z):=\delta_c$ be the point mass on some random string $c\in \{0,1\}^n$.
Let $P(X|Z)$ and $P(Y|Z)$ both be  given by the stochastic  matrix $A\otimes A \otimes \cdots  \otimes A$. Let $P_0\neq P_1$ be the probability vectors given by the columns of $A$.
Then 
\[
P(X)=P(Y)={\bf  P}_{c}\,,
\]
with ${\bf P}_c$ as in Definition~\ref{StringM}.
Since $P(Z)$ is supported by the singleton set $\{c\}$, we have $P(X|Y)=P(X)$ and $P(Y|X)=P(Y)$. Thus
\begin{eqnarray*}
K(P(X))+K(P(Y|X))&\Ceq &K(P(X|Y))+K(P(Y))\\&\Ceq &K(P(X))+K(P(Y))\Ceq 2 n\,.
\end{eqnarray*}
On the other hand, we have
\[
K(P(X|Z))+K(P(Y|Z))+K(P(Z))\Ceq 0+0+n=n\,.
\]

By observing that there is a third  variable $Z$ such that
\[
K(P(X|Z))+K(P(Y|Z))+K(P(Z))\Ceq K(P(X,Y))\,,
\]
we  thus have obtained a hint that the latent model is the more appropriate causal hypothesis.

\subsection*{Analysis of the required sample size}

The following arguments show that the above complexities of the Markov kernels become relevant already for moderate sample size.
Readers who are not interested in technical details may skip the remaining part of the subsection.

Consider first the sampling required to estimate $c$ by drawing i.i.d.~from ${\bf P}_c$ as in Definition~\ref{StringM}.
By counting the number of symbols $1$  that occur at position $j$ we  can guess whether
$c_j$ is $0$ or $1$ by choosing the distribution for which the relative frequency is closer  to the  corresponding
probability. 
To  bound the error probabilities from above
set 
\[
\mu:=|P_0(1)-P_1(1)|\,.
\]
Then the probability $q$ that the relative frequency deviates 
by more than $\mu/2$ decreases exponentially in the number  of  copies, i.e.,  
 $q \leq e^{-\mu m \alpha}$ where  $\alpha$ is an appropriate constant. 
The probability to have no error for any digit 
is then bounded from below by $(1-e^{-\mu  m  \alpha})^n$. We want to increase  $m$ such that
the error probability tends to zero. To this end, choose $m$ such that $e^{-\mu m \alpha} \leq 1/n^2$, i.e.,
$m \geq \ln n^2/(\mu \alpha)$. Hence 
\[
\Big(1-e^{-\mu  m  \alpha}\Big)^n \geq \Big(1-\frac{1}{n^2}\Big)^n \rightarrow 1
\]
because 
\[
\Big(1-1/n^2\Big)^{n^2} \rightarrow 1/e\,,
\]
and therefore
\[
\lim_{n \to \infty} \Big(1-1/n^2\Big)^n = \sqrt[n]{ \lim_{n\to\infty} \Big(1-1/n^2\Big)^{n^2}}=\lim_{n \to \infty} \sqrt[n]{1/e}=1\,.
\]
The required sample size thus grows only logarithmically in $n$.
In the same way, one shows that the sample size needed to distinguish between  
different conditionals $P(Y|X)={\bf A}_c$ increases only with the logarithm of $n$ provided that $P(X)$ is a strictly positive
product distribution on $\{0,1\}^n$.

\subsection{Resolving  statistical ensembles into individual observations}

\label{Subsec:Res}

The assumption of independent identically distributed random variables is one of the cornerstones of standard statistical reasoning.
In this section we show that the {\it independence} assumption in a typical
statistical sample
  is often due to 
prior knowledge on causal relations among single objects which can
nicely represented by a DAG. We will see that the algorithmic causal Markov condition then leads to non-trivial implications. 

\begin{figure}
\centerline{\includegraphics[scale=0.3]{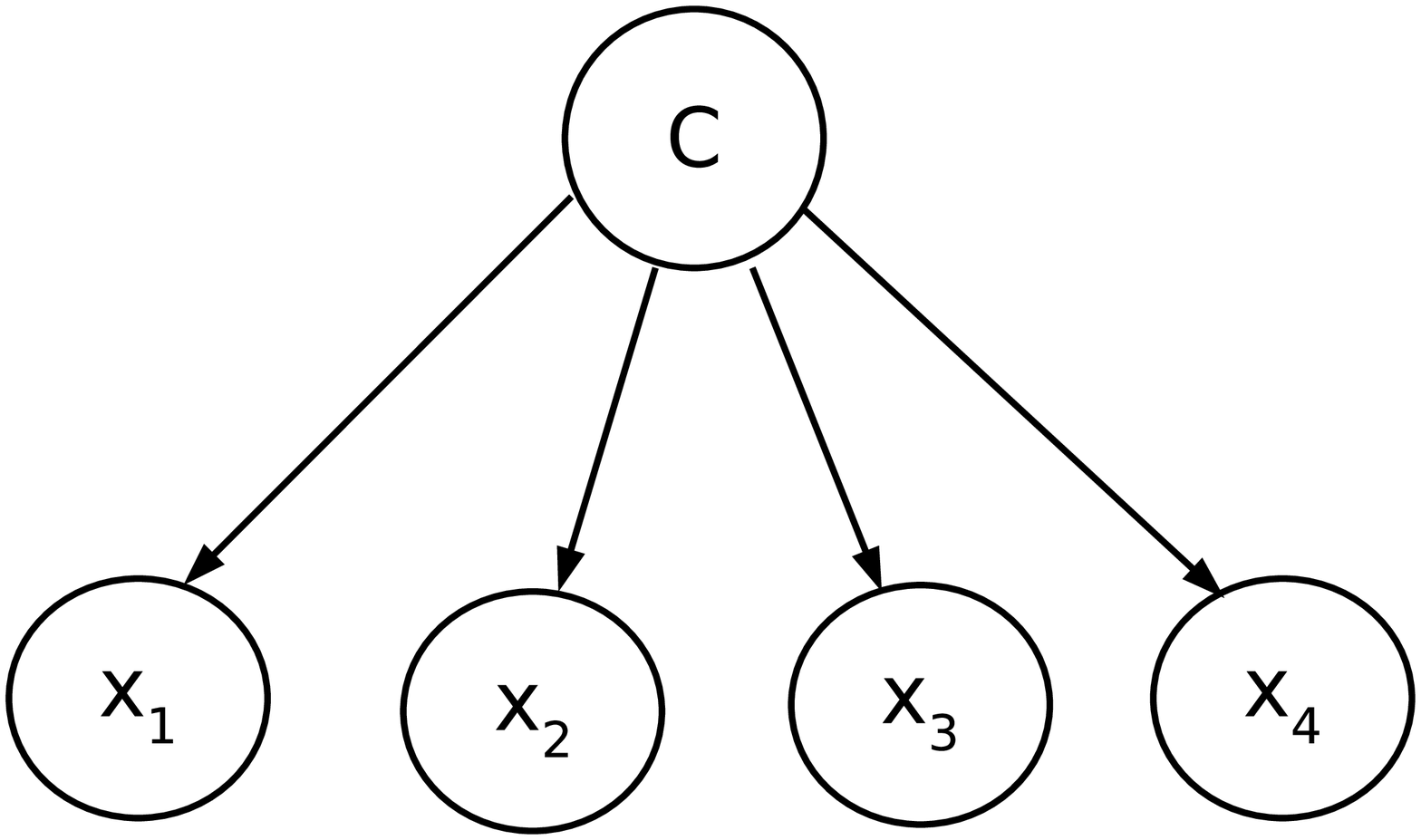}}
\caption{{\small Causal structure of the coin toss. The statistical properties of the coin $C$ define the common cause
that links the results of the coin toss.}\label{coin}} 
\end{figure}

 Assume we describe a biased coin toss, $m$ times repeated, 
and obtain the binary string $x_1,\dots,x_m$ as result. This  is certainly one of the scenarios where the i.i.d. assumption is well justified
because if we do not believe that the coin changes 
or that the result of one coin toss influences the other ones. The only relation between the coin tosses
is that they refer to  the same coin. We will thus draw  a DAG representing the relevant causal relations for the  scenario where
$C$ (the coin) is the common cause of all $x_j$ (see fig.~\ref{coin}).

Given the relevant information on $C$ (i.e., given the probability $p$ for ``head''),  we  have conditional algorithmic independence between the $x_j$ when
applying the Markov condition to this causal graph.\footnote{This is consistent 
with the following Bayesian interpretation: if we  define a non-trivial prior on the possible values of $p$, the 
individual observations are statistically dependent when marginalizing over the prior, 
but knowing $p$ renders  them
independent.}
However, there are two problems: (1)  
it does not make sense
to consider algorithmic mutual  information among binary strings of length $1$. (2) Our theory developed so far 
(Theorems~\ref{equiAMK} and \ref{algoImplM})
considered the  number  of strings
(which is $m+1$ here) as constant and
thus even the complexity of $x_1,\dots,x_m$ is considered as $O(1)$. 
To solve  this problem, we  define a new structure with three nodes as  follows. 
For some arbitrary $k<m$ set ${\bf x}^1:=x_1,\dots,x_k$ and ${\bf x}^2:=x_{k+1},\dots,x_m$.  
Then $C$ is the common cause of ${\bf x}^1$ and  ${\bf x}^2$ and
$I({\bf x}^1;{\bf x}^2|C)=0$ because every similarity between ${\bf x}^1$ and ${\bf x}^2$ is due to their common source
(note that
the information that the strings ${\bf x}^j$ have been obtained by combining $k$ and  $n-k$ results,
respectively, is here implicitly considered as background information in the sense of relative  causality in Subsection~\ref{Subsec:Rel}).
We will  later discuss examples where a source generates symbols from a larger  probability space. Then 
every $x_j$ is a string and 
it is important to keep in mind the ``format information'', i.e., the information how to read
the concatenation $x_1,  x_2, \cdots, x_k$ as a sample of $m$ strings. This format information will always be considered
as background, too. 

Of course, we may also consider partitions into more than two substrings keeping in mind that their number
is considered as  $O(1)$.
When we consider causal relations between {\it short} strings we will thus always 
apply the algorithmic causal Markov condition to groups of strings rather than applying it to the ``small objects''  itself. 
The DAG that formalizes the causal relations between
instances or groups of instances of a statistical ensemble and the source that determines the statistics
in the above sense
 will be  called the ``resolution of statistical  ensembles into individual observations''.

\begin{figure}
\centerline{\includegraphics[scale=0.3]{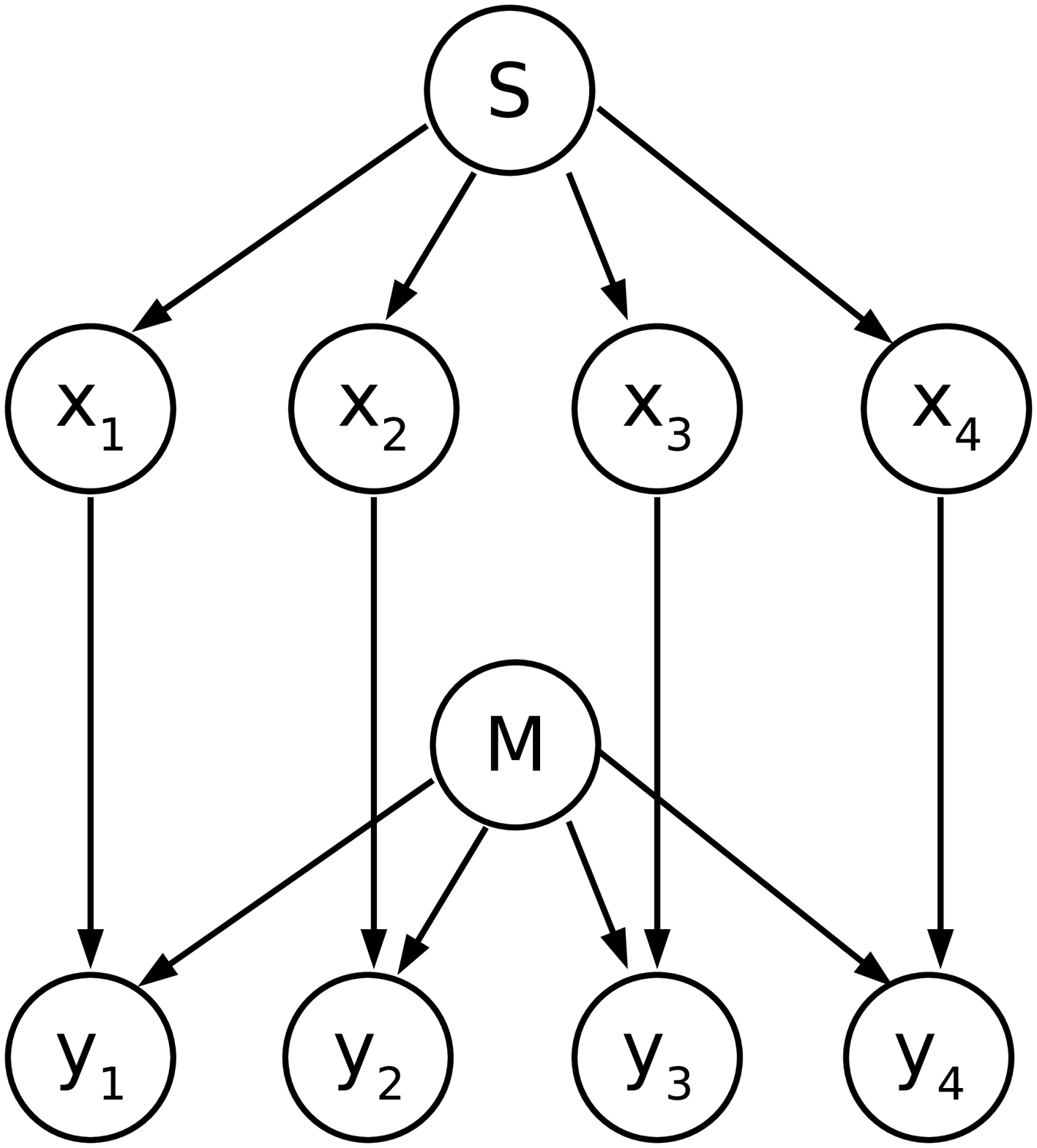}\hspace{1.0cm}\includegraphics[scale=0.3]{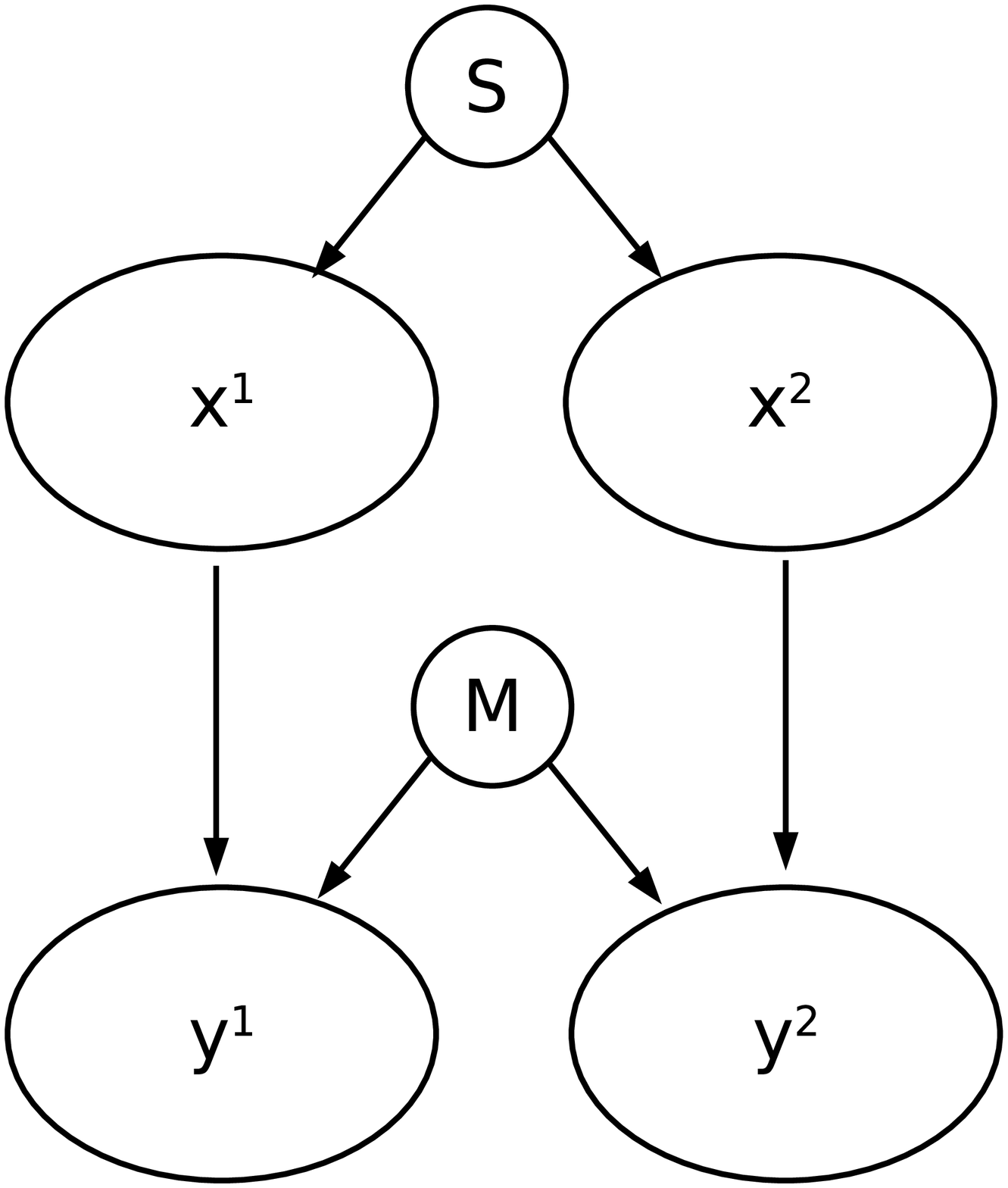}}
\caption{{\small Left:  causal structure obtained by resolving the causal structure $X\rightarrow Y$ 
between the random variables $X$ and $Y$ 
into causal relations among single events. Right: causal graph obtained by combining the first $k$ observations to ${\bf x}^1$ and the
remaining $m-k$ to ${\bf x}^2$ and  the same for $Y$. We observe that ${\bf x}^2$ d-separates ${\bf x}^1$ and ${\bf y}^2$, while ${\bf y}^2$ does not d-separate 
${\bf y}^1$ and ${\bf x}^2$. This asymmetry distinguishes causes from effects.}\label{SM}} 
\end{figure}

The resolution gets more interesting  if we consider causal relations between two random variables
  $X$ and $Y$. Consider the following  scenario where $X$ is  the cause of $Y$. Let $S$ be a source generating  $x$-values $x_1,\dots,x_m$ 
according to a fixed probability distribution $P(X)$. Let $M$ be a machine that receives these values as inputs and generates
$y$-values $y_1,\dots,y_m$ according to the conditional  $P(Y|X)$. Fig~\ref{SM} (left) shows  the causal graph  for $m=4$.

In analogy to the  procedure above,
we divide the string ${\bf x}:=x_1,\dots,x_m$ into ${\bf x}^1:=x_1,\dots,x_k$ and  ${\bf x}^2:=x_{k+1},\dots,x_m$ and  
use the same grouping for the $y$-values. 
We then draw the causal graph in fig.~\ref{SM}  (right) showing causal relations between ${\bf x}^1,{\bf x}^2,{\bf y}^1,{\bf y}^2,S,M$. Now we assume that $P(X)$ and $P(Y|X)$ 
are not known, i.e., we don't have access to  the relevant properties of $S$ and  $M$.  Thus we have to consider $S$ and $M$ as ``hidden
objects'' (in analogy to hidden variables  in the statistical setting). Therefore we have to apply the  Markov condition 
to the causal  structure  in such a way that only  the observed objects  ${\bf x}^1,{\bf x}^2,{\bf y}^1,{\bf y}^2$ occur. One checks easily
 that
${\bf x}^2$ d-separates ${\bf x}^1$ and ${\bf  y}^2$ and ${\bf x}^1$ d-separates ${\bf x}^2$ and ${\bf y}^1$.
Exhaustive search  over all possible triples of 
subsets of ${\bf x}^1,{\bf x}^2,{\bf y}^1,{\bf y}^2$ shows  
 that these are the only  non-trivial d-separation conditions. We conclude
\begin{equation}\label{asymcon}
I({\bf x}^1;{\bf y}^2|{\bf x}^2)\Ceq 0 \quad \hbox{ and }  \quad I({\bf x}^2;{\bf y}^1|{\bf x}^1)\Ceq 0    \,.
\end{equation}
The most remarkable property of eq.~(\ref{asymcon}) is that  it is asymmetric with respect to exchanging the roles of  $X$ and $Y$
since, for instance, $I({\bf y}^1;{\bf x}^2|{\bf y}^2)\Ceq 0$ can be violated. 
Intuitively, the reason is that given ${\bf y}^2$, the knowledge of ${\bf x}^2$ 
provides better insights into the  properties of $S$ and $M$ 
than knowledge of ${\bf x}^1$ would do,
which can be an advantage when describing ${\bf y}^1$.  
The following example shows  that this asymmetry can even be  relevant for sample size $m=2$ provided that the  probability space is 
large.

Let $S$ be a source that always generates the same string  $a\in \{0,1\}^n$. 
Assume  furthermore   that $a$ is algorithmically  random in  the sense that $K(a)\Ceq n$. 
For sample size $m=2$ we then have
${\bf x}=(x_1,x_2)=(a,a)$. Let $M$ be a machine that randomly removes $\ell$ digits
randomly either at the beginning or the end from its input string of length $n$. By this procedure we obtain 
a string $y_j\in  \{0,1\}^{\tilde{n}}$ with $\tilde{n}:=n-\ell$ from $x_j$.
  
\begin{figure}
\centerline{\includegraphics[scale=0.3]{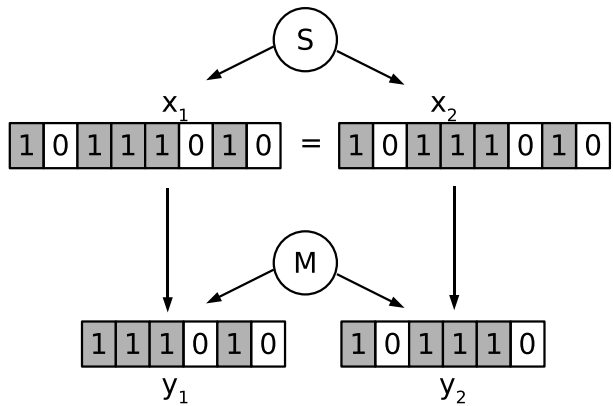}}
\caption{{\small Visualization of the truncation process: The  source  $S$  generates  always  the same  string, the machine truncates
either the left or the  right end. Given only the four strings $x_{1},x_2$ and $y_1,y_2$ as observations, we can reject the causal hypothesis
$Y\rightarrow X$. This is because $I(x_1:y_2|x_1)$ can be significantly greater than zero provided that the substrings missing in $y_1,y_2$ at
the left or at the right end, respectively, are sufficiently complex. \label{TwoStrings}}} 
\end{figure}

For sample size $2$ it is likely that 
$y_1$ and $y_2$ contain the last $n-\ell$ and the first $n-\ell$ digits of $a$, respectively, or vice versa. This process
is depicted in fig.~\ref{TwoStrings} for $n=8$ and  $\ell=2$.
Since the sample size is only two, the partition of the sample into two halves leads to single observations, i.e.,
${\bf x}^j=x_j$ and ${\bf y}^j=y_j$ for $j=1,2$.

 In short-hand notation, 
${\bf y}^1=a_{[1..n-\ell]}$ and ${\bf y}^2=a_{[\ell+1..n]}$. We then have
\[
I({\bf x}^1;{\bf y}^2|{\bf x}^2)\Ceq 0 \quad \hbox{ and } \quad I({\bf y}^1;{\bf x}^2|{\bf x}^1) \Ceq 0\,,
\]
but  
\[
I({\bf y}^1;{\bf x}^2|{\bf y}^2)\Ceq \ell \quad  \hbox{ and } \quad I({\bf x}^1;{\bf y}^2|{\bf y}^1)\Ceq \ell\,,
\]
which correctly lets us prefer the causal direction $X\rightarrow Y$ because these  dependences violate the global algorithmic Markov condition in Theorem~\ref{equiAMK}
when applied to a hypothetical  graph where ${\bf y}^1$ and ${\bf y}^2$ are the  outputs of the  source and ${\bf x}^1$ and ${\bf x}^2$ 
are the outputs of a machine that has received ${\bf y}^1$ and ${\bf y}^2$.

Even though the condition in eq.~(\ref{asymcon}) does not explicitly contain the notion of complexities of Markov kernels it is closely related to 
the algorithmic independence of Markov kernels. To explain this, assume we would generate algorithmic dependences between 
$S$ and $M$ by adding an arrow $S\rightarrow M$ or $S\leftarrow M$ or by adding a common cause. 
Then ${\bf x}^2$ would no longer d-separate ${\bf x}^1$ from ${\bf y}^2$. The possible violation of eq.~(\ref{asymcon})
could then be an observable result of the algorithmic dependences between the hidden objects $S$ and $M$ (and their statistical properties $P(X)$ 
and $P(Y|X)$, respectively).

\subsection{Conditional density estimation on subsamples}

\label{Subsec:MDL}

Now we develop an inference rule that is even closer to the idea of checking algorithmic dependences of Markov kernels
than condition~(\ref{asymcon}),
but still avoids the  notion of {\it Kolmogorov complexity of the ``true'' conditional distributions} by using
finite sample estimates instead.
Before we explain the idea we mention two simpler approaches for  doing so and describe  their potential problems.
It 
would be  straightforward  to  apply Postulate~\ref{inMech}
to the
finite sample estimates
of the conditionals.
In particular,
minimum description length (MDL) 
approaches
\cite{Gruenwald} appear promising from  the  theoretical point of view due to their close relation to Kolmogorov  complexity. 
We rephrase the minimum complexity estimator described by  Barron and Cover \cite{BarronCover}: 
Given a string-valued random variable $X$ and a sample $x_1,\dots,x_m$ 
drawn  from $P(X)$, set
\[
\hat{P}_m(X):= {\rm argmin} \Big\{ K(Q)- \sum_{j=1}^m \log Q(x_j)  \Big\}\,,
\] 
where $Q$ runs over all probability densities on the probability space  under consideration. 
If the data is sampled from a  computable  distribution, then
 $\hat{P}_m(X)$ converges in probability to $P(X)$ \cite{BarronCover}. Let us define a similar estimator  $\hat{P}_m(Y|X)$ for  the conditional density $P(Y|X)$. Could we reject the causal hypothesis $X\rightarrow Y$ after observing 
that $\hat{P}_m(X)$ and $\hat{P}_m(Y|X)$ are mutually  dependent? In the context of the true probabilities, we have argued that
$P(X)$ and  $P(Y|X)$ represent independent mechanisms. However, for the estimators we do not see a justification for 
independence because the relative  frequencies of the $x$-values influence the estimation of $\hat{P}_m(X)$ {\it and}
$\hat{P}_m(Y|X)$. 
This counter-argument becomes irrelevant only if the sample size is  such that the complexities of the estimators
coincide with the complexities of the true distributions. 
If we assume that the latter are  typically uncomputable (because generic real numbers are uncomputable) this sample size will never  be attained. 

The general idea of MDL  \cite{Gruenwald} 
also suggests the following causal inference  principle: 
If we are given the data points $(x_j,y_j)$ with $j=1,\dots,m$,
consider 
the MDL estimators $\hat{P}_m(X)$ and $\hat{P}_m(Y|X)$. They define a joint distribution that we denote by 
$\hat{P}_{X\rightarrow Y}(X,Y)$ (where we have  dropped $m$ for convenience).
The total description length 
\[
C_{X\rightarrow Y}:= K(\hat{P}_m(X))+K(\hat{P}_m(Y|X)) - \sum_{j=1}^m\log\hat{P}_{X\rightarrow Y}(x_j,y_j)
\]
measures the complexity of the probabilistic model plus the complexity of the data,  given  the model.
Then we  compare $C_{X\rightarrow Y}$ to $C_{Y\rightarrow X}$ (defined  correspondingly) and prefer  
the causal direction with the smaller value. However, it is not clear whether this kind of reasoning can  be  derived
from the algorithmic Markov condition.

For this reason, we construct an inference rule that uses estimators in a more sophisticated way and whose
justification is directly based on applying  the algorithmic Markov condition to the resolution of ensembles.
The idea of our strategy is that we do not use the full data set to estimate $P(Y|X)$. 
Instead, we apply the estimator to a subsample of $(x,y)$ pairs
that no longer carries significant information about the relative frequencies of $x$-values in the full data set. As we  will see below,
this  leads to
 algorithmically independent {\it finite sample} estimators for the Markov kernels if  the causal hypothesis is correct.

\begin{figure}
\centerline{\includegraphics[scale=0.38]{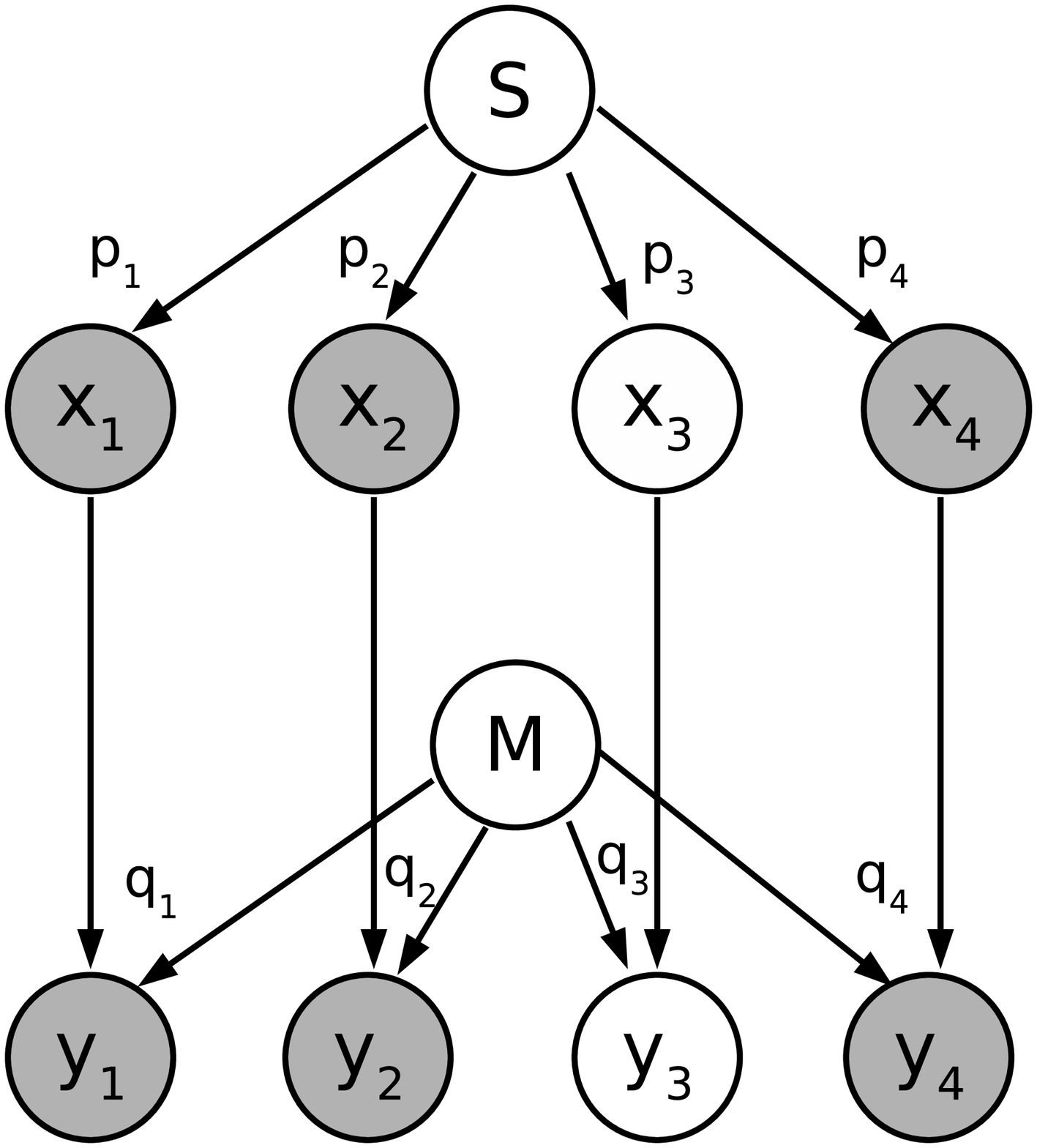}\hspace{1cm}\includegraphics[scale=0.38]{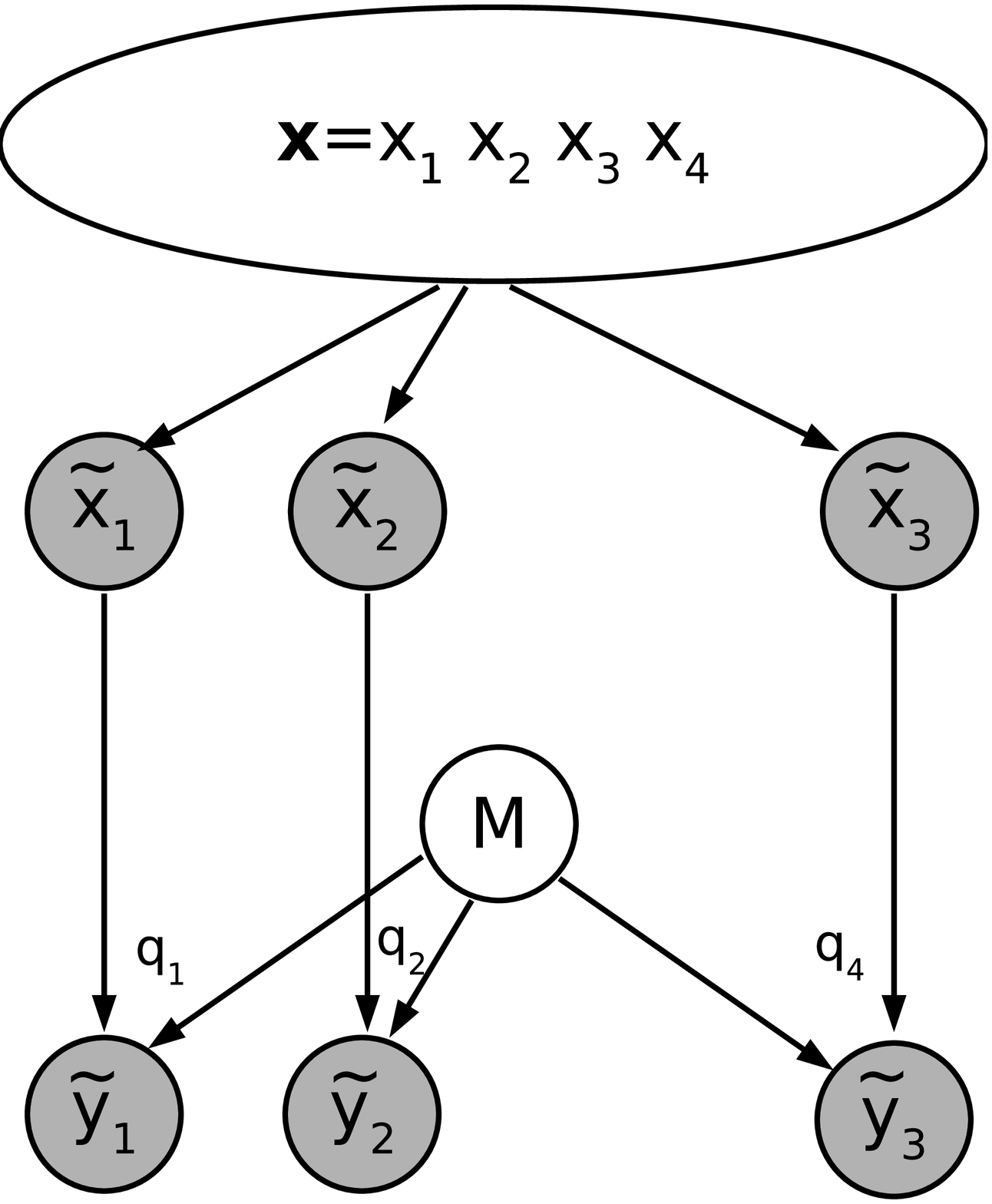}}
\caption{{\small (left) Causal structure between single observations $x_1,\dots,x_m$,$y_1,\dots,y_m$
for sampling  from  $P(X,Y)$, given the causal structure $X\rightarrow  Y$. 
The  programs $p_j$ compute $x_j$ from the description of the source $S$. The programs  $q_j$  compute $y_j$ from $x_j$  and the 
description of the machine $M$, respectively. 
The  grey nodes are those that are selected  for the subsample (see text). Right:  
Causal structure relating ${\bf x}$, $\tilde{x}_j$, and $\tilde{y}_j$.
Note that the causal relation  between $\tilde{x}_j$ and $\tilde{y}_j$ is the same as  the one between
the corresponding  pair $x_j$  and $y_j$. Here, for instance,  $\tilde{x}_3=x_4$ and  $\tilde{y}_3=y_4$ 
and it is thus still the same program $q_4$ that computes 
$y_4$ from $x_4$ and  $M$. 
Hence, the causal model that links $M$ with the selected values $\tilde{x}_j$ and $\tilde{y}_j$ is the subgraph
of the graph showing relations between $x_j$, $y_j$ and $M$. This kind of robustness of the causal structure with respect to the  selection procedure will be used below.}\label{SM_q}} 
\end{figure}

Let $X\rightarrow Y$ be the causal structure that generated the 
data  $({\bf x},{\bf  y})$,
with ${\bf x}:=x_1,\dots,x_m$ and ${\bf y}:=y_1,\dots,y_m$  after $m$-fold i.i.d.~sampling from $P(X,Y)$.
The resolution of the ensemble is the causal graph in fig.~\ref{SM_q}, left. 

According to Postulate~\ref{algoFunc} there are mutually independent programs
$p_j$ computing $x_j$ from the description of $S$. Likewise, there are mutually independent  programs  
$q_j$ computing $y_j$ from $M$ and $x_j$.
Assume we are given a rule  how  to generate a subsample of $x_1,\dots,x_m$ from ${\bf x}$.
It is important  that this selection rule does not refer  to  ${\bf y}$ but only  uses ${\bf x}$ (as well as some
random string as additional  input)
and that 
the selection can be performed by a  program of length $O(1)$. 
Denote the subsample by
\[
{\bf \tilde{x}}=\tilde{x}_1,\dots,\tilde{x}_l:=x_{j_1},\dots,x_{j_l}\,,
\]
with $l<m$.  
The  above selection of indices defines also a subsample of $y$-values
\[
{\bf y}:=y_{j_1},\dots,y_{j_l}:=\tilde{y}_1,\dots,\tilde{y}_l\,.
\] 
By construction, we  have
\[
\tilde{y}_i=p_{j_i}(\tilde{x}_i,M)\,.
\]
Hence we can draw the causal structure depicted in fig.~\ref{SM_q}, right.

\begin{figure}
\centerline{\includegraphics[scale=0.38]{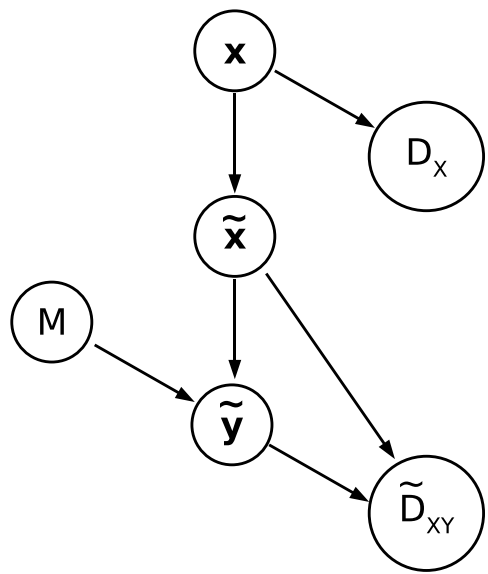}}
\caption{{\small $D_X$  is some information  derived from ${\bf x}$. 
The idea is  that it is a density estimator for $P(X)$ or that it describes properties  of  the empirical distribution of $x$-values.
If
the selection procedure  ${\bf x} \rightarrow \tilde{{\bf x}}$ has sufficiently blurred this information,
the mutual information between $\tilde{ {\bf x}}$ and  $D_X$ is  low.  
$D_{XY}$ on the other hand, is a density estimator for $P(Y|X)$ or it 
encodes some desired properties of the empirical joint distribution of $x$- and $y$-values 
in the subsample.
If the mutual information between $D_X$ and $\tilde{D}_{X,Y}$ exceeds the one between $\tilde{{\bf x}}$ and  $D_X$, 
we reject the hypothesis $X\rightarrow Y$.}\label{SM_D}} 
\end{figure}

Let now $D_X$ be any string that is derived from ${\bf  x}$ by some program of length  $O(1)$. 
$D_X$ 
may be the full description of relative frequencies or any {\it computable} density estimator $\hat{P}(X)$,
or some other  description of interesting properties of the relative frequencies. 
Similarly,  
let $\tilde{D}_{YX}$ be a description that is derived from ${\bf x},{\bf y}$ by some simple algorithmic rule.
The idea is that it is a computable estimator $\tilde{P}(Y|X)$ for the conditional distribution $P(Y|X)$ or any relevant property of the latter.  Instead of estimating conditionals, one may also consider an estimator of the  {\it joint} density of the subsample.
We augment  the causal structure in fig.~\ref{SM_q}, right, with $D_X$ and  $\tilde{D}_{YX}$.
The structure can be simplified  by  merging nodes in the same  level and we obtain
the structure in fig.~\ref{SM_D}.

To derive testable implications of the causal hypothesis,
we  observe  that every information between  $D_X$ and $\tilde{D}_{YX}$ is processed via ${\bf \tilde{x}}$. We thus have
\begin{equation}\label{indD}
\tilde{D}_{XY} \independent D_X \, | {\bf \tilde{x}}^* \,,
\end{equation}
which formally follows  from the global Markov condition in Theorem~\ref{equiAMK}.
Using Lemma~\ref{dataPr} and eq.~(\ref{indD})  we conclude
\begin{equation}\label{DataProcessD}
I(D_X;\tilde{D}_{YX})\stackrel{+}{\leq} I({\bf \tilde{x}}; D_X)\,.
\end{equation}

The intention behind generating the subsample ${\bf \tilde{x}}$ is to ``blur'' the distribution of  $X$.
If we have a density estimator $\hat{P}(X)$ we try to choose  the subsample such that the algorithmic mutual information between
${\bf \tilde{x}}$ and $\hat{P}(X)$ is small. Otherwise we have not sufficiently blurred the distribution of $X$.
Then we apply an arbitrary conditional density estimator $\hat{P}(Y|X)$ to the subsample. If  there still is 
a non-negligible amount of mutual information
between $\hat{P}_X$ and $\hat{P}(Y|X)$, the causal hypothesis in fig.~\ref{SM}, left, cannot be true
and we reject $X\rightarrow Y$.

To show that the above procedure can also be applied to data sampled from {\it uncomputable} probability distributions, 
let  $P_0$ and $P_1$ be uncomputable distributions on $\{0,1\}$
and $A_0,A_1$ uncomputable stochastic maps  from $\{0,1\}$ to $\{0,1\}$.
Define a  string-valued random variable $X$ with distribution $P(X):={\bf P}_c$ 
as in Definition~\ref{StringM} and the conditional distribution of a string-valued variable $Y$ by
$P(Y|X):={\bf A}_d$ as in Definition~\ref{StringJ} for  strings $c,d\in \{0,1\}^n$. 
Let $P_0$ and $P_1$ as well as $A_0$ and $A_1$ be
known up to an accuracy that is sufficient to distinguish between them. 
We assume that all this information (including $n$)  is given as background knowledge, but $c$ and $d$ are unknown.
Let  $D_X=:c'$, where $c'$ is the estimated value of $c$ computed from the finite sample ${\bf x}$ of size $m$.
Likewise, let $\tilde{D}_{XY}:=d'$ be the estimated value of $d$ derived from the subsample $({\bf x},{\bf y})$ of size $\tilde{m}$.
If $m$ is large enough (such that also $\tilde{m}$  is sufficiently large) we can estimate $c$ and $d$, i.e,  $c'=c$  and $d'=d$  with high probability. 
The most radical method  to blur $P(X)$ is to choose
 ${\bf \tilde{x}}$ such that the empirical distribution of $x$-values is  uniform and
the $x_j$-values  are lexicographically reordered (with some random  ordering among the $j$-values that  correspond to the same $x$-value).    
The only algorithmic information that ${\bf \tilde{x}}$ then contains is the description of
its  length, i.e., $\log_2 \tilde{m}$ bits. Hence we have
\[
I(D_X:\tilde{{\bf x}})\stackrel{+}{\leq} \log_2 \tilde{m}\,.
\] 
Assume now that $c=d$. Then 
\[
I(D_X:\tilde{D}_{XY})\Ceq n\,,
\]
provided that the estimation was correct.
As shown at the end of Subsection~\ref{Subsec:plMKAlg},  
this is already possible for  $\tilde{m}=O(\log n)$, i.e.,
\[
I(D_X:{\bf\tilde{x}})\in O(\log_2 n)\,,
\] 
which violates ineq.~(\ref{DataProcessD}). The importance of  this example lies in the fact that $I(P(X):P(Y|X))$  is not well-defined here
because $P(X)$ and $P(Y|X)$ both are uncomputable. Nevertheless, $P(X)$ and $P(Y|X)$ have a computable aspect, i.e, 
the strings $c$ and $d$ characterizing them. Our strategy is therefore suitable to detect algorithmic dependences between 
computable features.

It is remarkable that the above scheme 
is general enough to include  also strategies for very small sample sizes provided  that the
probability space is large. 
To describe an extreme case, we  consider again the example with the truncated strings in 
fig.~\ref{TwoStrings} with the role of $X$ and  $Y$ 
reversed.  Let $Y$ be a random variable whose value is always the constant string $a\in \{0,1\}^n$.
Let $P(X|Y)$ be the mechanism that generates $X$ by truncating 
either the $l$ leftmost digits or the $l$ rightmost digits of $Y$ (each with probability $1/2$).
We denote these  strings by $a_{{\rm left}}$ and $a_{{\rm right}}$, respectively.
Assume  we have two observations $x_1=a_{{\rm left}}$, $y_1=c$ and $x_2=a_{{\rm right}}$, $y_2=a$.
We define a subsample by selecting 
only the first observation $\tilde{x}_1:=x_1=a_{{\rm left}}$ and $\tilde{y}_1:=y_1=a$. Then we define $D_X:=x_1, x_2$ 
and $\tilde{D}_{XY}:=y_1$. We observe that the mutual information between $\tilde{D}_{XY}$ and $D_X$ is $K(a)$, 
while the mutual information between $D_X$ and ${\bf \tilde{x}}$ is only $K(a_{{\rm left}})$.
Given generic choices of $a$, 
this violates condition  (\ref{DataProcessD})  and we reject the causal hypothesis $X\rightarrow Y$.

\subsection{Plausible Markov  kernels in  time series}

\label{Subsec:TS}

Time series are  interesting examples of causal structures   
where the time order provides prior knowledge on the causal direction. 
Since there is a large number of them available from all scientific disciplines they can be useful to test causal inference rules on data with known 
ground truth. 
Let us consider  the following  example of a
causal inference problem. 
Given a  time series  and the prior knowledge that  it has been generated by a first order Markov process, but the
direction is unknown. Formally, we are given observations $x_1,x_2,x_3,\dots,x_m$ corresponding to random variables  $X_1,X_2,\dots,X_m$ 
such that the causal structure is either
\begin{equation}\label{true}
\cdots \rightarrow X_1\rightarrow X_2 \rightarrow X_3 \cdots \rightarrow X_n \rightarrow \cdots \,,
\end{equation}
or
\begin{equation}\label{false}
\cdots \leftarrow  X_1 \leftarrow  X_2\leftarrow X_3 \cdots \leftarrow X_n \leftarrow \cdots \,,
\end{equation}
where we have extended  the series to infinity in both directions. 

The question is whether the asymmetry of the joint distribution with respect to time  inversion provides hints
on the real time direction. 
Let us assume now that the graph (\ref{true}) corresponds to the true time direction. Then 
the  hope is that $P(X_{j+1}|X_j)$ is simpler, in some reasonable sense,  than $P(X_j|X_{j+1})$. 
At first glance this seems to be a straightforward extension of the principle of plausible Markov kernel discussed in Subsection~\ref{Subsec:plMKAlg}.
However, there is a subtlety with the justification when we apply our ideas to stationary time series:

Recall that the principle of minimizing the total complexity of all  Markov kernels over all potential causal directions
has been derived from the independence of  the true Markov kernels 
(remarks after Postulate~\ref{inMech}). 
However,  the algorithmic independence of 
 $P(X_j|PA_j)=P(X_j|X_{j-1})$ and $P(X_i |PA_i)=P(X_i|X_{i-1})$ fails spectacularly 
because
 stationarity implies that these Markov kernels {\it  coincide} and represent a causal mechanism that is constant in time. 
This shows that  the justification of minimizing  total complexity breaks down for stationary   time series. 

The following argument shows that not only the justification breaks down but also the principle as such:
Consider the case where $P(X_j)$ is the  unique stationary distribution of the Markov  kernel
$P(X_{j+1}|X_j)$. Then we have
\[
K(P(X_j|X_{j+1})) \stackrel{+}{\leq} K(P(X_{j+1},X_j))\Ceq K(P(X_{j+1}|X_j))\,.
\]
Because the forward time conditional describes uniquely the backward  time   conditional   (via implying the
description of the unique stationary marginal) the Kolmogorov complexity of the latter can  exceed the complexity of the former only by a constant term.

We now focus on {\it non}-stationary time series. To motivate the general idea we first present  an example described in \cite{OccamsRazor}.
Consider a random walk of a particle on $\Z$ starting at  $z\in \Z$. In  every  time step the  probability is $q$ to move one 
site to the right and  $(1-q)$ to move to  the left. Let $X_j$ with  $j=0,1,\dots$ 
be the random  variable describing the position after step $j$. Then we have $P(X_0=z)=1$.
The forward time conditional reads 
\[
P(x_{j+1}|x_j)= \left\{ \begin{array}{ccc} q & \hbox{ for } &  x_{j+1}=x_j+1 \\
 1-q & \hbox{ for } & x_{j+1}=x_j-1 \\
0 & \hbox{  otherwise \,.} & \end{array} \right. 
\]  
To   compute the backward time conditional we first compute $P(X_j)$ which is given by the distribution of a Bernoulli experiment with $j$ steps.
Let $k$ denote  the number  of right moves, i.e., $j-k$ is the number of left moves. 
With 
$x_j=k-(j-k)+z=2k-j+z$ we thus obtain 
\[
P(x_j)=  q^{(j+x_j-z)/2} (1-q)^{(j-x_j+z)/2} { j \choose (j+x_j-z)/2}\,. 
\]  
Elementary calculations show
\begin{eqnarray*}
P(x_j|x_{j+1})&=&P(x_{j+1}|x_j)\frac{P(x_j)}{P(x_{j+1})}\\ &=& 
\left\{\begin{array}{ccc}   \frac{(j+x_j-z)/2 +1}{j+1}          & \hbox{ for }  &  
x_j=x_{j+1}-1 \\  
  \frac{(j-x_j+z)/2+1}{j+1}   & \hbox{ for }  & x_j=x_{j=1}+1  \\ 0 &  &\hbox{ otherwise\,.}   \end{array}\right.
\end{eqnarray*}
The forward time  process is specified  by  the initial condition
$P(X_0)$ (given by $z$) and the transition probabilities $P(X_j,\dots,X_1|X_0)$ 
(given by $p$). 
A priori, these two ``objects'' are mutually unrelated, i.e.,
\begin{eqnarray*}
&&K(P(X_0),P(X_j,X_{j-1},\dots,X_1|X_0))\Ceq \\
&&K(P(X_0))+K(P(X_j,X_{j-1},\dots,X_1|X_0))\Ceq \\&&K(z)+K(q)\,.
\end{eqnarray*}
On the other hand, 
the description of $P(X_j)$ (the ``initial condition'' of the backward time process) 
alone already requires the specification of {\it both}  $z$ and $q$. The
description of the  ``transition rule'' $P(X_1,\dots,X_{j-1}|X_j)$ refers only to $z$. 
We thus have
\[
K(P(X_j))+K(P(X_1,X_{2},\dots,X_{j-1}|X_j))\Ceq 2 K(z)+K(q)\,.
\]
Hence 
\[
I(P(X_j):P(X_0,X_1,\dots,X_{j-1}|X_j)) \Ceq K(z)\,.
\]
The fact that the initial distribution of the hypothetical process 
\[
X_j \rightarrow X_{j-1} \rightarrow \cdots \rightarrow X_0
\] 
shares algorithmic information with the transition probabilities makes the hypothesis suspicious.

\subsection*{Resolving time series}

We have seen that the algorithmic dependence between ``initial condition'' and ``transition rule'' of the backward time process
(which would be surprising if it occurred for the forward time process)
represents an asymmetry of non-stationary time-series with respect  to time reflection.
We will now discuss this asymmetry after resolving the statistical  ensemble into individual observations.

\begin{figure}
\centerline{\includegraphics[scale=0.4]{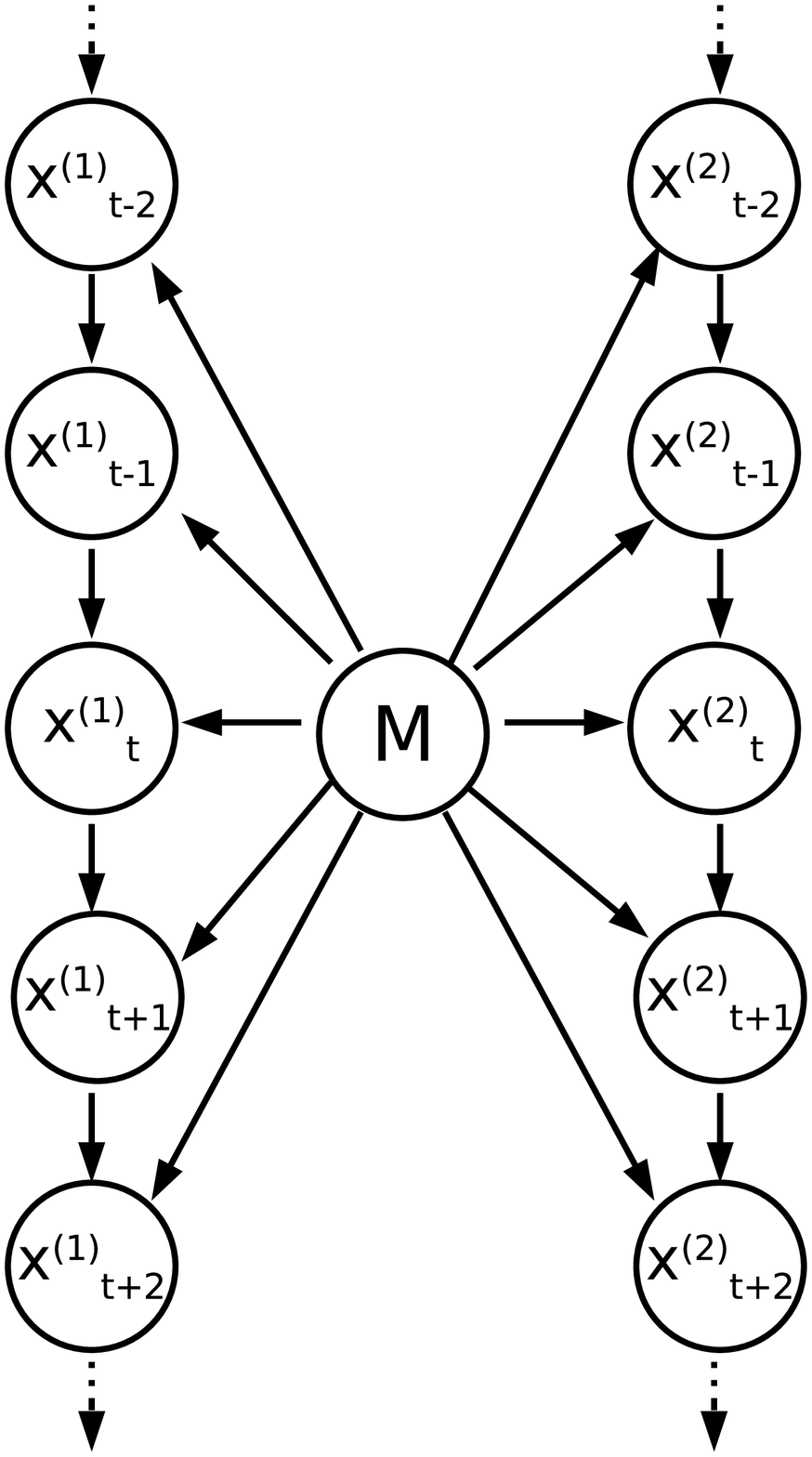}\hspace{0.5cm}  \includegraphics[scale=0.4]{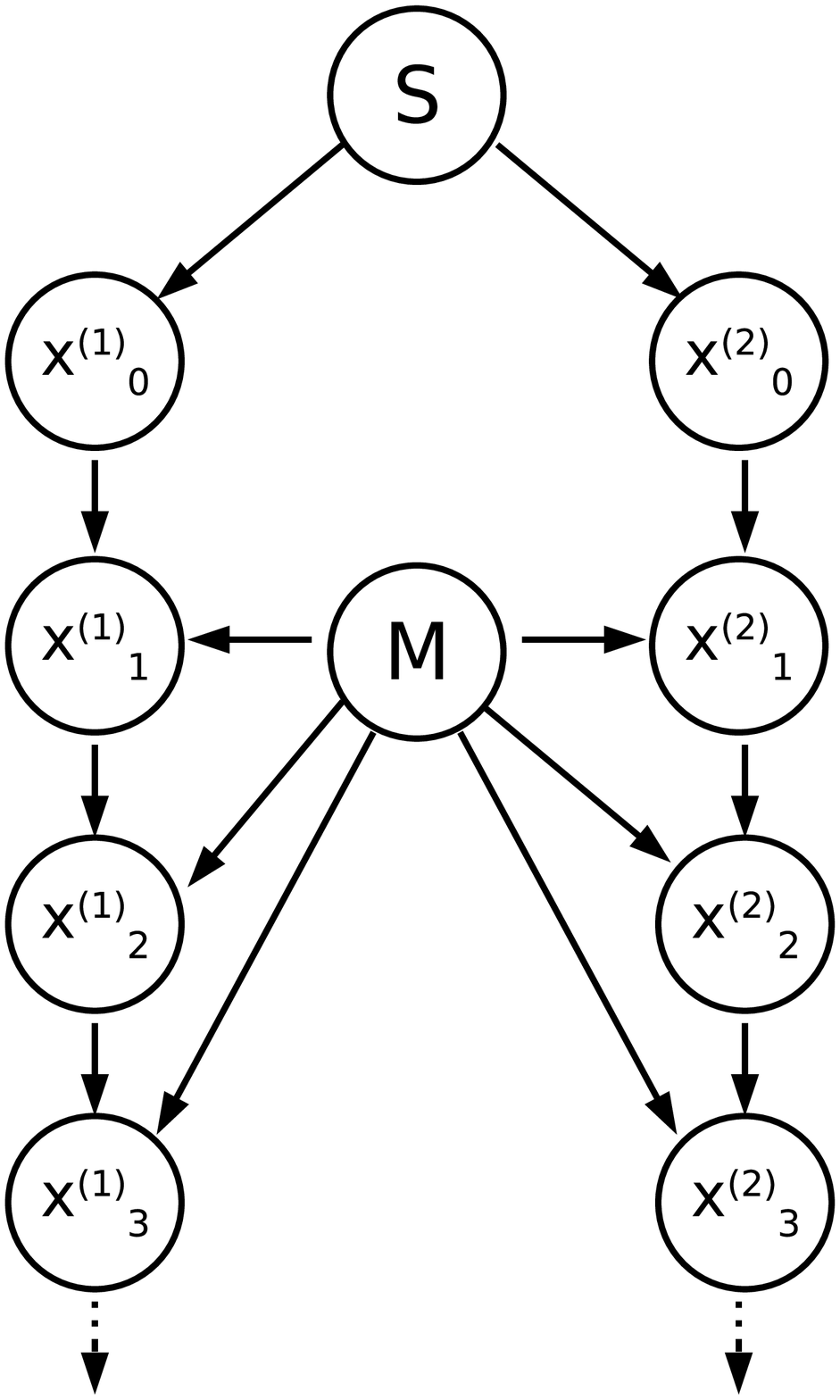} }
\caption{{\small Left: causal graph of a time series. The values $x^{(j)}_i$ corresponds  to the $j$th instance at  time $i$.
Right: the initial part of the time-series is asymmetric with respect to time-inversion.  
\label{timeseries}}} 
\end{figure}

Assume we are given $m$ instances of $n$-tuples $x^{(i)}_1,\dots,x^{(i)}_n$  with $i=1,\dots,m$ that have been  
i.i.d.~sampled from $P(X_1,\dots,X_n)$ and $X_1,\dots,X_n$ are part of a time series that can be 
described by a  first order stationary Markov  process.
Our resolution of a statistical ensemble generated by $X\rightarrow Y$ contained a source $S$ and a machine $M$. The  source 
generates $x$-values
and the machine generates $y$-values from the input  $x$. The algorithmic independence of $S$ and $M$ was essential for the 
asymmetry between cause and effect described in Subsection~\ref{Subsec:Res}.
For the  causal chain 
\[
\cdots \rightarrow X_1\rightarrow X_2 \rightarrow X_3 \rightarrow \cdots
\]
we would therefore have machines $M_j$ generating the $x_j$-value from $x_{j-1}$. However, for 
stationary time-series all $M_j$ are the {\it same} machine. 
The causal structure of the resolution of the statistical ensemble for $m=2$ is shown in fig.~\ref{timeseries}, left.

This graph entails no independence constraint that is asymmetric with  respect 
to reversing the time direction. To see this, recall that two DAGs entail the same set of
independences if and only if they have the same skeleton (i.e. the corresponding undirected graphs coincide) and the same set of unshielded colliders ($v$-structures), i.e., substructures  
$A \rightarrow C \leftarrow B$ where $A$ and $B$ are non-adjacent \cite{Pearl:00}.
Fig.~\ref{timeseries} has no such $v$-structure and the skeleton is obviously symmetric with respect  to time-inversion.

The initial part  is, however, asymmetric (in agreement with the asymmetries entailed by fig.~\ref{SM}, left) 
and we have 
\[
I(x^{(1)}_0:x^{(2)}_1|x^{(2)}_1)\Ceq 0\,.
\]
This is just the finite-sample analogue of the  statement that the initial distribution $P(X_0)$ and the transition rule 
$P(X_j|X_{j-1})$ are algorithmically
independent.

\section{Decidable modifications of the inference rule}

\label{Sec:Dec}

To use the algorithmic Markov condition in practical applications
we have to replace it with {\it computable} notions  of complexity.
The following  two subsections discuss two different directions along  which practical inference 
rules can be developed.

\subsection{Causal inference using symmetry constraints}

\label{Subsec:Sym}

We have seen that the algorithmic causal Markov condition implies that the the sum of the  
Kolmogorov complexities of the Markov kernels must be minimized over all possible causal graphs. 
In practical applications, it is natural to replace  the minimization of Kolmogorov complexity with a decidable simplicity criterion
even though this makes the relation to the theory developed so far rather vague. In this subsection we will describe
an empirically decidable inference rule and show that the relation to Kolmogorov complexity of conditionals
is closer  than it may seem  at first glance. 

Moreover, the example below shows a scenario where the  causal  hypothesis $X\rightarrow Y$ can already be
preferred to
 $Y\rightarrow X$ by comparing only the {\it marginal} distributions $P(X)$ and $P(Y)$ and observing that
a simple conditional $P(Y|X)$ leads from the former to the  latter but no {\it simple} conditional leads  into the
opposite direction.
The example will furthermore show
why  the  identification of causal directions is often easier for {\it probabilistic} causal relations than for 
{\it  deterministic} ones, a point that has  also been pointed  
out by Pearl \cite{Pearl:00} in a different context.

\begin{figure}
\centerline{\includegraphics[scale=0.45]{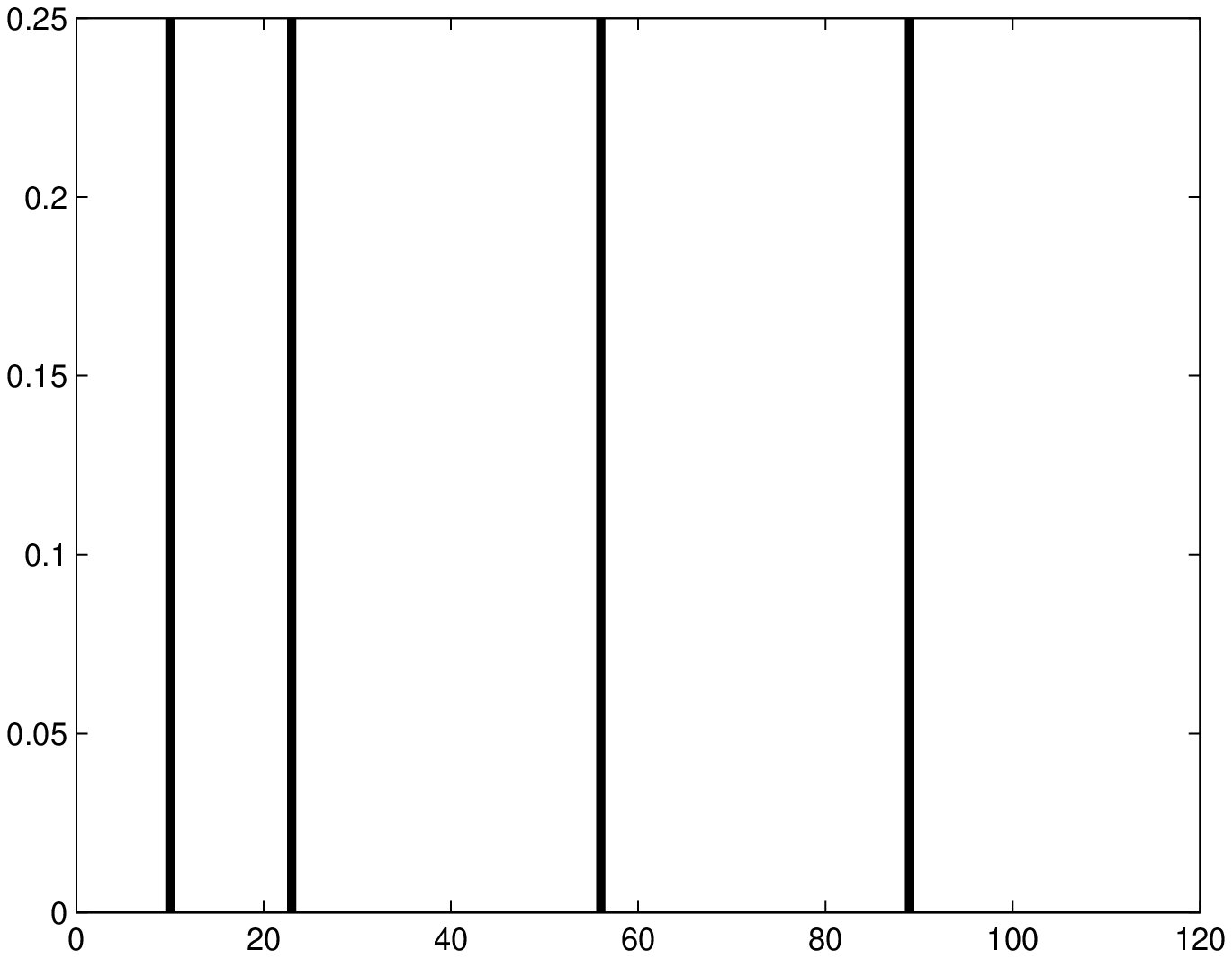}\includegraphics[scale=0.45]{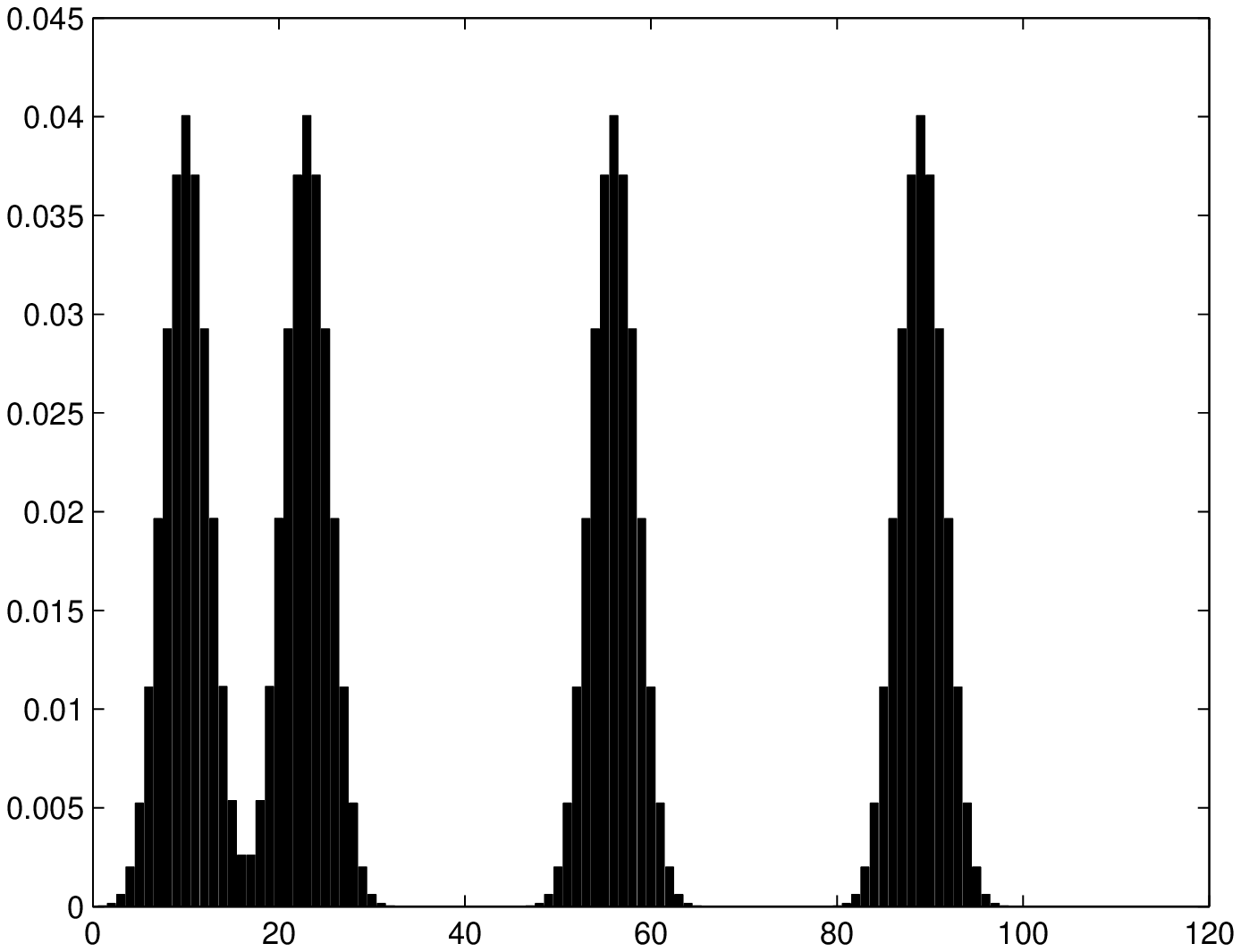}}
\caption{\label{peaks}{\small Two probability distributions $P(X)$ (left) and $P(Y)$ (right) 
on the set $\{1,\dots,120\}$ both having $4$ peaks at the positions $n_1,\dots,n_4$,
but the peaks in $P(X)$ are well-localized and those of  $P(Y)$ are smeared out by a  random walk}}
\end{figure}

Consider  the discrete probability space  
$\{1,\dots,N\}$. Given two distributions $P(X),P(Y)$  like the  ones depicted in
fig.~\ref{peaks} for $N=120$. The   marginal $P(X)$ consists of $k$ sharp peaks of equal height  at positions $n_1,\dots,n_k$
and  $P(Y)$ also has $k$ modes centered at the same positions, but with greater width.
We assume that $P(Y)$ can be  obtained from $P(X)$ by repeatedly applying a doubly stochastic matrix $A=(a_{ij})_{i,j=1,\dots,N}$
with $a_{ii}=1-2p$ for $p\in (0,1)$ and $a_{ij}=p$ for $i=j \pm 1 ({\rm mod}\, N) $. 
The stochastic map $A$ thus defines a random walk and we have by assumption 
\[
P(Y)=A^m P(X)
\]
for some $m\in\N$. 

Now  we ask which causal hypothesis is more likely: (1) $P(Y)$ has been obtained from $P(X)$ 
by some stochastic map $M$.
(2) $P(X)$ has been obtained from   $P(Y)$  by some stochastic map $\tilde{M}$.
Our  assumptions already contain  an example $M$ that corresponds to  the first hypothesis ($M:=A^m$).
Clearly, there also  exist  maps $\tilde{M}$ for hypothesis (2). One example would be  
\begin{equation}\label{StrMtilde}
\tilde{M}:=[P(X),P(X),\dots,P(X)]\,,
\end{equation} 
i.e. $M$
has the probability vector $P(X)$ in every column. 

To describe in  which sense $X\rightarrow Y$  is the simpler hypothesis we observe that
$\tilde{M}$ in eq.~(\ref{StrMtilde}) already contains the description of the 
positions $n_1,\dots,n_k$  whereas
$M=A^m$ is rather simple. 
The Kolmogorov complexity of $\tilde{M}$  as chosen  above is
for a generic choice of  the positions $n_1,\dots,n_k$ given  by 
\[
K(\tilde{M})\Ceq K(P(Y))\Ceq \log {N \choose k}\,,
\]
where $\Ceq$ denotes equality up to a term that does not depend on $N$.
This is because 
different locations $n_1,\dots,n_k$ of  the original peaks lead to different distributions  $P(Y)$ and,
conversely, every such $P(Y)$ is uniquely 
defined by describing the positions of the corresponding sharp peaks and $M$.

However, we want to prove that also other choices of $\tilde{M}$ necessarily have high values of Kolmogorov  complexity. 
To this end, we define a  family of ${N \choose k}$
probability distributions $P_j(X)$ given by equally high peaks at the positions $n_1,\dots,n_k$ and accordingly
the smoothed probability distributions $P_j(Y)$. We first need the  following result.

\begin{Lemma}[average complexity of stochastic maps]${}$\\
\label{data_proc}
Let $(Q_j(X))_{j=1,\dots,\ell}$ and $(Q_j(Y))_{j=1,\dots,\ell}$ be 
two families of marginal distributions of $X$ and $Y$, respectively.
 
Moreover, let $(A_j)_{j=1,\dots\ell}$ be  a family of not necessarily different stochastic matrices
with $A_jQ_j(Y)=Q_j(X)$.
Then 
\[
\frac{1}{\ell} \sum_{j=1}^\ell K(A_j)  \geq I(X:J)-I(Y:J)\,,
\]
where
the information that $X$ contains about the index $j$ is given by
\[
I(X:J):=H\Big(\frac{1}{\ell}\sum_j Q_j(X)\Big) -\frac{1}{\ell}\sum_j H(Q_j(X))\,,
\]
$J$ denotes the random variable with values $j$. Here, $H(.)$ denotes the Shannon entropy and  $I(Y:J)$ is computed in a similar way as
$I(X:J)$ using $Q_j(Y)$ instead  of $Q_j(X)$.
\end{Lemma}

\vspace{0.3cm}
\noindent  
Proof: The idea is to show that we need at least $2^{\Delta}$ different stochastic matrices to achieve that
the information $I(X:J)$ exceeds $I(Y:J)$ by the amount $\Delta$. Using a standard argument 
rephrased below,
the average complexity is therefore at  least $\Delta$. 

Assume, for instance, that all $A_j$ coincide. Then the usual data processing inequality \cite{Cover} shows that 
applying the same matrix  to the different distributions 
can never increase the information on the index $j$, i.e., $I(X:J)\leq I(Y:J)$.  
To derive the lower bound on the number of {\it different} matrices required
we  define a partition of $\{1,\dots,\ell\}$ into  $d$ sets
$S_1,\dots,S_k$ for which the $A_j$ coincide.
In other words, we have
 $A_j=B_r$ if $j\in S_r$ and the matrices $B_1,\dots,B_d$ are chosen appropriately. 
We define a random variable $R$ whose value $r$ indicates that $j$ lies in the $r$th equivalence class.
The above ``data processing argument'' implies
\begin{equation}\label{sameE}
I(X:J|R)\leq I(Y:J|R)\,.
\end{equation}
Furthermore, we have 
\begin{equation}\label{differentE}
I(X:R)\leq I(Y:R)+\log_2 d\,.
\end{equation}
This is because both $I(X:R)$ and $I(Y:R)$ cannot exceed $\log_2 d$ because $d$ is the number of  values $R$ can attain. 
Then we have:
\begin{eqnarray*}
I(X:J)&=&I(X:J,R)=I(X:R)+I(X:J|R)\\&\leq& I(Y:R)+\log_2 d+I(Y:J|R)\\&=& \log_2 d+I(Y:J)\,.
\end{eqnarray*} 
The first equality follows because $R$ contains no {\it additional} information on $X$ 
(when $J$ is known) since it describes only from which equivalence class $j$ is taken.  
The second equality  is a general rule for mutual information \cite{Cover}. 
The inequality combines ineqs.~(\ref{sameE}) and (\ref{differentE}).
The last equality follows similar as the equalities in the first line. This shows that we need at least
$2^d$ different matrices with $d=\lceil I(X:J)-I(Y:J) \rceil $.
We have
\[
2^{-\frac{1}{d} \sum_{j=1}^d  K(B_j)} \leq \frac{1}{d} \sum_{j=1}^d 2^{-K(B_j)} \leq \frac{1}{d}\,,
\]
where the first inequality holds because the exponential function is concave  and the second is entailed by Kraft's inequality.
This yields
\[
\frac{1}{d} \sum_{j=1}^d  K(B_j)  \geq \log_2 d\,,
\]
completing the proof.
$\Box$

\vspace{0.3cm}
\noindent
To
apply  Lemma~\ref{data_proc} to the above example we define families of $\ell:={N \choose k}$ 
distributions  $P_j(X)$ having their peaks  at the positions
$n_1,\dots,n_k$ and also their smoothed versions $P_j(Y)$.
Mixing all probability distributions will generate the entropy $\log N$ for $P_j(X)$ because we then obtain the uniform distribution.
Since we  have assumed that $P_j(Y)$ can be obtained from $P_j(X)$ by a doubly stochastic map, 
mixing
 all $P_j(Y)$ also  yields the  uniform distribution.
Hence the difference between  $I(X:J)$ and $I(Y:J)$ is simply given by the average entropy difference 
\[
\Delta H:=\frac{1}{\ell}  \sum_{j=1}^\ell \Big( H(P_j(Y))- H(P_j(X))\Big)\,.
\]
The Kolmogorov complexity required to map $P_j(Y)$ to $P_j(X)$ is thus, on average  over all $j$, at least
the entropy generated by the double stochastic random walk. Hence we  have shown that a typical example of two distributions
with peaks at  arbitrary positions $n_1,\dots,n_k$ needs a process $\tilde{M}$ whose Kolmogorov complexity 
is at least the entropy difference. 

One may ask why  to  consider  distributions with several peaks even though  the above result
will formally also apply to distributions $P_j(X)$ and $P_j(Y)$ with only {\it one} peak. 
The problem is  that the statement ``two distributions 
have a peak at the same  position'' does  not necessarily make sense for empirical data. This is because the definition of variables
is often chosen such that  the distribution becomes centralized. The statement that {\it multiple} peaks occur on seemingly
random positions seems therefore more sensible than the statement that {\it one} peak has been observed at a random position. 

\begin{figure}
\centerline{\includegraphics[scale=0.6]{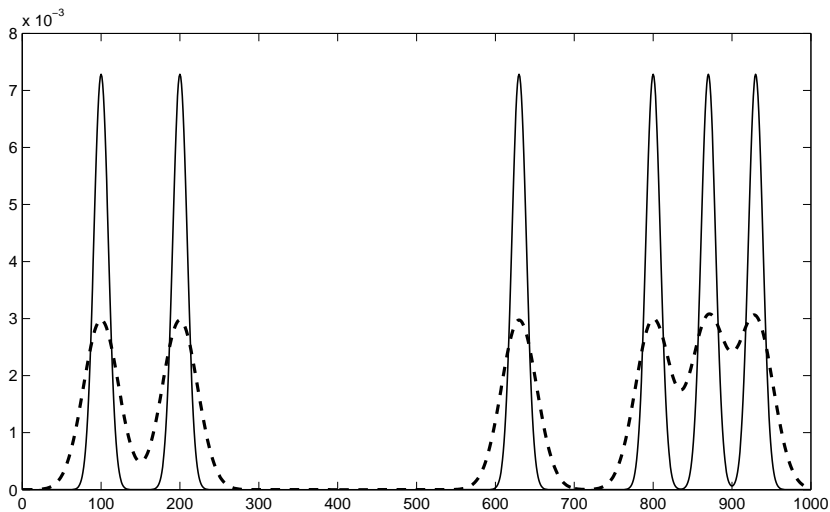}}
\caption{\label{smooth_peaks}{\small Two probability distributions $P(X)$ (solid) and $P(Y)$ (dashed) 
where $P(Y)$  can be  obtained from $P(X)$ by convolution with
a Gaussian  distribution}}
\end{figure}

We have above used  a finite number of discrete bins in order or keep the problem  as much combinatorial as possible.
In reality, we would rather expect a scenario like  the one in fig.~\ref{smooth_peaks} where two distributions
on $\R$ have  the same  peaks,  but the peaks  in the one distribution have been smoothed, for example by an additive Gaussian noise.

As above, we would rather assume that $X$ is the cause of $Y$ than vice versa since the  smoothing process  is simpler
than any process  that leads in the opposite direction. We emphasize that {\it denoising} is  an operation  that cannot be  
represented by a stochastic matrix, it is a linear operation that can be applied to the whole data set 
in order to reconstruct the original peaks. 
The statement is thus that no simple {\it stochastic process} leads  in the opposite direction. 
To  further discuss the rationale behind this way of reasoning we introduce another notion of
simplicity that  does not refer to  Kolmogorov complexity.
To this end, we  introduce the  notion of translation covariant 
conditional probabilities:

\begin{Definition}[translation covariance]${}$\\ 
Let  $X,Y$ be two real-valued random variables. 
A conditional distribution $P(Y|X)$ with density $P(y|x)$ is called translation covariant if
\[
P(y|x+t)=P(y-t|x)\quad \quad \forall t\in \R\,.
\]
\end{Definition}

Apart from this, we will also need the following well-known concept from 
statistical 
estimation theory \cite{Cra}:

\begin{Definition}[Fisher information]${}$\\
Let $p(x)$ be a continuously differentiable 
probability density of $P(X)$ on $\R$. Then the Fisher information 
is defined as
\[
F(P(X)):=\int \left( \frac{d}{dx} \ln p(x) \right)^2 dx\,.
\]
\end{Definition}

Then we have the following  Lemma (see 
Lemma~1  in  \cite{clock} showing the statement in a more general
setting that  involves also quantum stochastic maps):

\begin{Lemma}[monotonicity under covariant maps]${}$\\
Let   $P(X,Y)$ be a joint  distribution such that $P(Y|X)$ is
translation covariant. Then
\[
F(P(Y))\leq F(P(X))\,.
\] 
\end{Lemma}

The intuition is that $F$ quantifies the  degree to which a  distribution is
non-invariant with respect to translations 
and that no translation covariant process
is able  to increase  this measure. The convolution with a Gaussian
distribution with non-zero variance decreases the 
Fisher  information. Hence there  is never a translation  invariant 
stochastic map in backward direction. 

The argument above  can easily be generalized in two respects. First, the argument works
also with other  quantities that  are monotonous with respect to 
translation invariant stochastic  maps.  Second, we  can also consider 
more general symmetries:

\begin{Definition}[general group covariance]${}$\\
Let $X,Y$ be random variables with equal range $S$.
Let $G$ be a group of bijections 
$g:S\rightarrow S$ and $X^g$ and $Y^g$ denoting
the random variables obtained by permuting the outcomes of the corresponding
random experiment according to $g$. 
Then we call a conditional $P(Y|X)$ $G$-covariant if
\[
P(Y^g|X))=P(Y|X^{g^{-1}})\quad \quad \forall g \in G\,.
\]
\end{Definition}

It is easy  to see that covariant stochastic maps define a  quasi-order of
probability distributions on $S$ by defining  $P\geq \tilde{P}$ if
there is a covariant stochastic map $A$ such that $A*P=\tilde{P}$. This 
is transitive  since the concatenation of covariant maps
is again covariant. 

If a $G$-invariant measure $\mu$ (``Haar measure'') exists on $G$ we can easily define an information
theoretic quantity that
measures the degree of non-invariance  with respect to $G$:

\begin{Definition}[reference information]${}$\\
Let $P(X)$ be a distribution on  $S$ and $G$ be  a group of bijections on $S$ with Haar measure $\mu$. 
Then the reference information is given by:
\begin{eqnarray*}
I_G&:=& H\left(P\Big[\int_G  X^g d\mu(g)\Big]\right)- \int_G  H\Big(P(X^g)\Big) d\mu(g)\\&=&
H\left(P\Big[\int_G  X^g d\mu(g)\Big]\right)- H(P(X))\,.
\end{eqnarray*}
\end{Definition}

The name ``reference information'' has been used
in \cite{GroupCovariantThermo} in a slightly different context where this  
information occurred as the value of a physical system to communicate 
a reference system (e.g. spatial or  temporal) where $G$
describes,  for instance,  translations in time or space.    
The quantity
$I_G$ can easily  be interpreted as mutual information $I(X:Z)$  if we 
introduce a $G$-valued random variable $Z$ whose values 
indicate which transformation $g$ has been applied. 
One can thus show that $I_G$ is non-increasing with  respect to every 
$G$-covariant map \cite{GroupCovariantThermo,Referenz}.

The following model describes a link between inferring causal directions by preferring covariant conditionals
to preferring directions with algorithmically independent Markov kernels.
Consider first the probability space  $S:=\{0,1\}$. We define the  
group $G:=\Z_2=(\{0,1\}, \oplus)$, i.e.,  the additive group of integers modulo $2$, acting
on $S$ as bit-flips or identity.
For any distribution on $P$  on $\{0,1\}$, the reference information $I_G(P)$ then measures the asymmetry with respect
to bit-flips. For two distributions $P$ and $\tilde{P}$ we can have the situation that a $G$-symmetric stochastic matrix
leads from $P$ to $\tilde{P}$, but only asymmetric stochastic maps convert  $\tilde{P}$ into $P$.  
Now we extend this idea to the group $\Z_2^n$ acting on strings of length $n$ by independent bit-flips. 
Assume we have a distribution on $\{0,1\}^n$  of the  form ${\bf P}_c$ in Definition~\ref{StringM} for some string $c$
and generate the distribution ${\bf \tilde{P}}_c$ by applying $M$ to ${\bf P}_c$  where
\[
M:=\left(\begin{array}{cc} 1-\epsilon_1 &\epsilon_1 \\ \epsilon_1 & 1-\epsilon_1 \end{array}\right) \otimes 
\left(\begin{array}{cc} 1-\epsilon_2 &\epsilon_2 \\ \epsilon_2 & 1-\epsilon_2 \end{array}\right)\otimes 
\cdots \otimes  \left(\begin{array}{cc} 1-\epsilon_n &\epsilon_n \\ \epsilon_n & 1-\epsilon_n \end{array}\right)\,,
\]
with  $\epsilon_j\in  (0,1)$. Then $M$ is  $G$-symmetric,  but
no $G$-symmetric process leads  backwards. This is because 
every such stochastic map would be asymmetric in a way that encodes $c$, i.e., the map
would have ``to know''  $c$ because $M$ has destroyed some amount of information about it.

\subsection{Resource-bounded  complexity}

\label{Subsec:RB}

The problem that the presence  or absence of mutual information
is undecidable
(when defined
 via Kolmogorov complexities) 
is similar to statistics, but  also 
different in other  respects. Let us first focus on the analogy.  
Given two real-valued random variables
$X,Y$, it is impossible to show by finite sampling that they are statistically independent.
$X \independent Y$ is equivalent to $E(f(X)g(Y))=E(f(X))E(g(Y))$ for every pair $(f,g)$ of measurable functions.
If  we observe significant correlations between $f(X)$ and  $g(Y)$ for some previously defined
pair, it  is well-justified to reject independence. The same holds if such correlations are detected for $f,g$ in some
previously defined,  sufficiently small set of functions (cf.~\cite{Gretton2005b}).  However, if this  is not  the case,
we can never be  sure that there is not some pair of
arbitrarily complex functions $f,g$ 
that are correlated  with respect to  the true distribution.
Likewise, if we have two strings   $x,y$ and find no simple program that computes $x$ from $y$  this does not  mean that there is no such a rule. Hence, we also have  the statement that there can  be an algorithmic dependence even though we do not find it.

However, the difference to the statistical situation is  the following. Given that we have found functions $f,g$ yielding correlations
it is only a matter  of the statistical significance level whether this is sufficient to reject independence. For algorithmic 
dependences, we do not even have a decidable criterion to reject {\it in}dependence. 
Given that we have found a simple  program that computes $x$ from  $y$,   it still may 
 be true that $I(x;y)$ is small because there may also be a simple rule to generate $x$  (which would imply $I(x:y)\approx 0$) that we were not able to find.  
This shows that we can neither show dependence nor independence. 

One possible answer to these problems is that Kolmogorov  complexity is only an idealization of empirically  decidable quantities. 
Developing this idealization only aims at  providing hints in which directions we have  to develop practical inference rules.
Compression algorithms have already been developed that are intended  to approximate, for  instance, 
the algorithmic information of genetic sequences \cite{grumbach94new,ChenCompression}. 
Chen et al.~\cite{ChenCompression}
constructed a ``conditional compression scheme'' to approximate conditional Kolmogorov complexity and 
applied it to the estimation of
the algorithmic mutual information between two genetic sequences. To evaluate to which extent  methods
of this kind can be  used for  causal inference using  the algorithmic Markov condition is an interesting 
subject  of  further research.

It  is also  noteworthy that  there is a theory on 
{\it resource-bounded} description complexity \cite{Vitanyi97} where compressions of $x$ are only allowed
if the decompression can be performed within a previously defined number of computation steps and on a tape of previously defined length.
An important advantage of resource-bounded complexity is  that it is computable. 
The disadvantage, on the other hand, is that the mathematical theory is more difficult. Parts of this paper have been developed
by converting statements on statistical dependences into their algorithmic counterpart. The strong analogy between statistical and
algorithmic mutual information occurs only for complexity with unbounded resources. For instance, the symmetry 
$I(x:y)\Ceq I(y:x)$ breaks  down when replacing  Kolmogorov complexity with resource-bounded versions \cite{Vitanyi97}. 
Nevertheless,
to develop a  theory of inferred causation using {\it resource-bounded} complexity could be a challenge for the future. 
There are several reasons to believe that taking into account computational complexity can provide additional hints on the causal structure:

Bennett \cite{BennettCom,BennettInt,BennettDepth}, for instance, 
has argued that the {\it logical depth} of an object echoes in some sense its history.
The former is, roughly speaking, defined as follows. Let $x$ be a string that describes the object
and $s$ be its shortest description. Then the logical depth of $x$ is the number of time  steps that a parallel
computing device requires  to compute $x$ from $s$. According to Bennett, large logical depth indicate that the object has been created by a 
process that consisted of many non-trivial steps. This  would mean  that there  also is some causal information that follows
from the time-resources required to compute a string from its shortest description. 

The time-resources required to compute one observation from the other also plays a role in
the discussion of causal inference rules in \cite{OccamsRazor}.
The paper presents a model where the conditional 
\[
P({\rm effect}|{\rm cause})
\]
can be  {\it efficiently} computed,  while
computing 
\[
P({\rm cause}|{\rm effect})
\]
 is NP-hard. This suggests that the computation time required to use information of the cause
for the description of the effect can be different from the time needed to obtain information on the cause from the effect. 
However, the goal of the  present paper was to  describe asymmetries between
cause  and  effect that even occur when computational complexity is ignored.

\section{Conclusions}

We have shown that 
our algorithmic causal Markov condition  links algorithmic dependences between single observations with the underlying causal structure in the same way
This is similar to the way the statistical causal Markov condition links statistical dependences among random variables to the causal structure.
The algorithmic Markov condition 
 has implications on different levels:

\vspace{0.2cm}
\noindent
(1) In conventional causal inference one can drop the assumption that observations  
\[
(x^{(i)}_1,\dots,x^{(i)}_n)
\]
 have 
been generated by {\it independent} sampling from a constant joint distribution 
\[
P(X_1,\dots,X_n)
\]
of  $n$ random variables $X_1,\dots,X_n$.
Algorithmic information theory thus replaces statistical causal inference with a probability-free formulation.

\vspace{0.2cm}
\noindent
(2) Causal relations among individual objects can be inferred provided their shortest descriptions are sufficiently complex.

\vspace{0.2cm}
\noindent
(3) New statistical causal inference rules follow  because causal hypotheses are suspicious if the corresponding Markov kernels
are algorithmically dependent.

\vspace{0.2cm}
Since algorithmic mutual information is uncomputable because Kolmogorov complexity is uncomputable,  we have
presented decidable inference rules that are motivated by the uncomputable idealization.



\end{document}